\documentclass[11pt, a4paper]{article}
\usepackage{amsmath}
\usepackage{amsfonts}
\usepackage{amssymb}
\usepackage{ulem}

\usepackage{amsthm}
\usepackage{url}

\usepackage{cite}
\usepackage{fancyhdr}

\usepackage{graphicx, color}
\usepackage{hyperref}

\newtheorem*{ques*}{Question}

\newtheorem*{rem*}{Remarks}

\begin{document}

\title{Convex Pentagons with Positive Heesch Number}
\author{ Teruhisa SUGIMOTO$^{ 1), 2)}$ }
\date{}
\maketitle

{\footnotesize

\begin{center}
$^{1)}$ The Interdisciplinary Institute of Science, Technology and Art

$^{2)}$ Japan Tessellation Design Association

E-mail: ismsugi@gmail.com
\end{center}

}

\medskip

{\small
\begin{abstract}
\noindent
We found convex pentagons whose Heesch number is equal to 
one, and which admit an edge-to-edge corona. In this manuscript, we present 
a new classification of these convex pentagons.
\end{abstract}

\textbf{Keywords: }convex pentagon, corona, 
Heesch number, monohedral tiling, tessellation

}

\section{Introduction}
\label{section1}

\begin{figure}[htbp]
 \centering\includegraphics[width=15cm,clip]{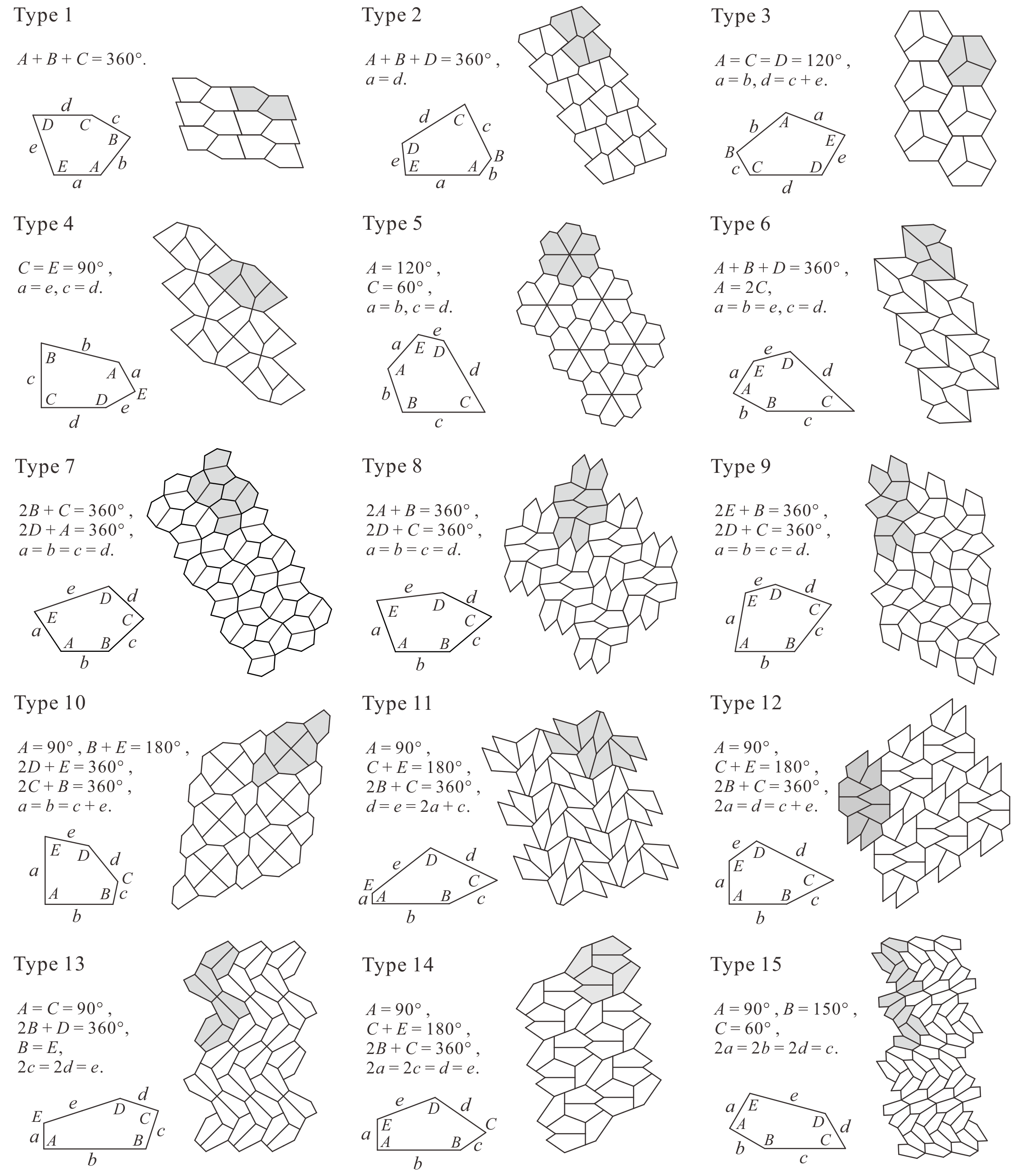} 
  \caption{{\small 
Convex pentagonal tiles of 15 families. Each of the convex pentagonal 
tiles is defined by some conditions between the lengths of the edges and the 
magnitudes of the angles, but some degrees of freedom remain. For example, a 
convex pentagonal tile belonging to Type 1 satisfies that the sum of three 
consecutive angles is equal to $360^ \circ$. This condition for Type 1 is 
expressed as $A+B+C=360^ \circ$ in this figure. The pentagonal tiles of Types 14 
and 15 have one degree of freedom, that of size. For example, the value of 
$C$ of the pentagonal tile of Type 14 is $\cos ^{ - 1}((3\sqrt {57} - 17) / 16) 
\approx 1.2099\;$rad $ \approx 69.32^ \circ $. The pale gray pentagons in 
each tiling indicate a fundamental region (the unit that can generate a 
periodic tiling by translation only).} 
\label{fig01}
}
\end{figure}

A planar tiling (or tessellation) is a collection of sets that are called 
tiles, which covers a plane without gaps and overlaps, except for the 
boundaries of the tiles. The term ``tile'' refers to a 
topological disk, whose boundary is a simple closed curve. If all tiles in 
the tiling are congruent, then the tiling is called \textit{monohedral}. Then, the 
polygon in the monohedral tiling is called the \textit{prototile} of monohedral 
tiling, or simply, the \textit{polygonal tile}~\cite{G_and_S_1987, Sugimoto_2012}. 
To date, fifteen families of convex pentagonal tiles, each of them 
referred to as a ``Type'' (i.e., Type 1, Type 2, etc. up to Type 15) are 
known\footnote{ In May 2017, Micha\"{e}l Rao declared that the complete 
list of Types of convex pentagonal tiles had been obtained (i.e., they have 
only the known 15 families), but it does not seem to be fixed as of January 
2018~\cite{Rao_2017, Wiki_PenTP}.} (see Figure~\ref{fig01}\footnote{ The 
classification of Types of convex pentagonal tiles is based on the essentially 
different properties of pentagons. The conditions of each Type express 
the essential properties. The classification problem of Types of convex 
pentagonal tiles and the classification problem of pentagonal tilings 
are quite different. The Types are not necessarily ``disjoint," that is, 
convex pentagonal tiles belonging to some Types also exist~\cite{Sugimoto_2012, 
Sugimoto_2016}.})~\cite{G_and_S_1987, Mann_2015, Rao_2017, Sugimoto_2012, 
Sugimoto_2017, Wiki_PenTP}. Moreover, a tiling by convex polygons is called 
edge-to-edge if any two convex polygons either are disjoint, or share 
one vertex or one entire edge. Within the representative tilings of the 
known families of convex pentagonal tiles, there exist edge-to-edge 
as well as non-edge-to-edge tilings~\cite{Sugimoto_2012, 
Sugimoto_2017, Wiki_PenTP}. For example, in Figure~\ref{fig01}, 
the representative tiling of Type 6 is edge-to-edge, whereas 
the representative tiling of Type 10 is non-edge-to-edge.
It is known that a convex pentagonal tile that can generate an 
edge-to-edge tiling, belongs to at least one of the 
known eight families~\cite{Bagina_2011, Sugimoto_2016}.

\renewcommand{\figurename}{{\small Figure.}}

Let $T$ be a tile in a plane. A \textit{corona} of $T$ is the set comprising of the centrally 
placed original $T$ and its layer of congruent copies of $T$. The corona is formed 
without gaps and overlaps (except for the boundaries of $T)$. The first corona 
is the set of all tiles sharing a boundary point with a centrally placed tile (including 
the original tile itself). For  an integer $h \ge 2$, the $h$-th corona is 
the set of all tiles sharing a boundary point with the $(h - 1)$-th corona 
(including the $(h - 1)$-th corona itself). The maximum number of layers, 
whose corona of $T$ can be formed, is called the \textit{Heesch number} of $T$ and is 
denoted by $H(T)$. If $T$ is a prototile of monohedral tiling, then $H(T) = \infty $. 
Thus, the convex pentagonal tiles are convex pentagons with $H(T) = \infty $. On 
the other hand, the regular pentagon is a convex pentagon with $H(T) = 0$ 
since it cannot generate a tiling and a corona.

\begin{figure}[htbp]
 \centering\includegraphics[width=13cm,clip]{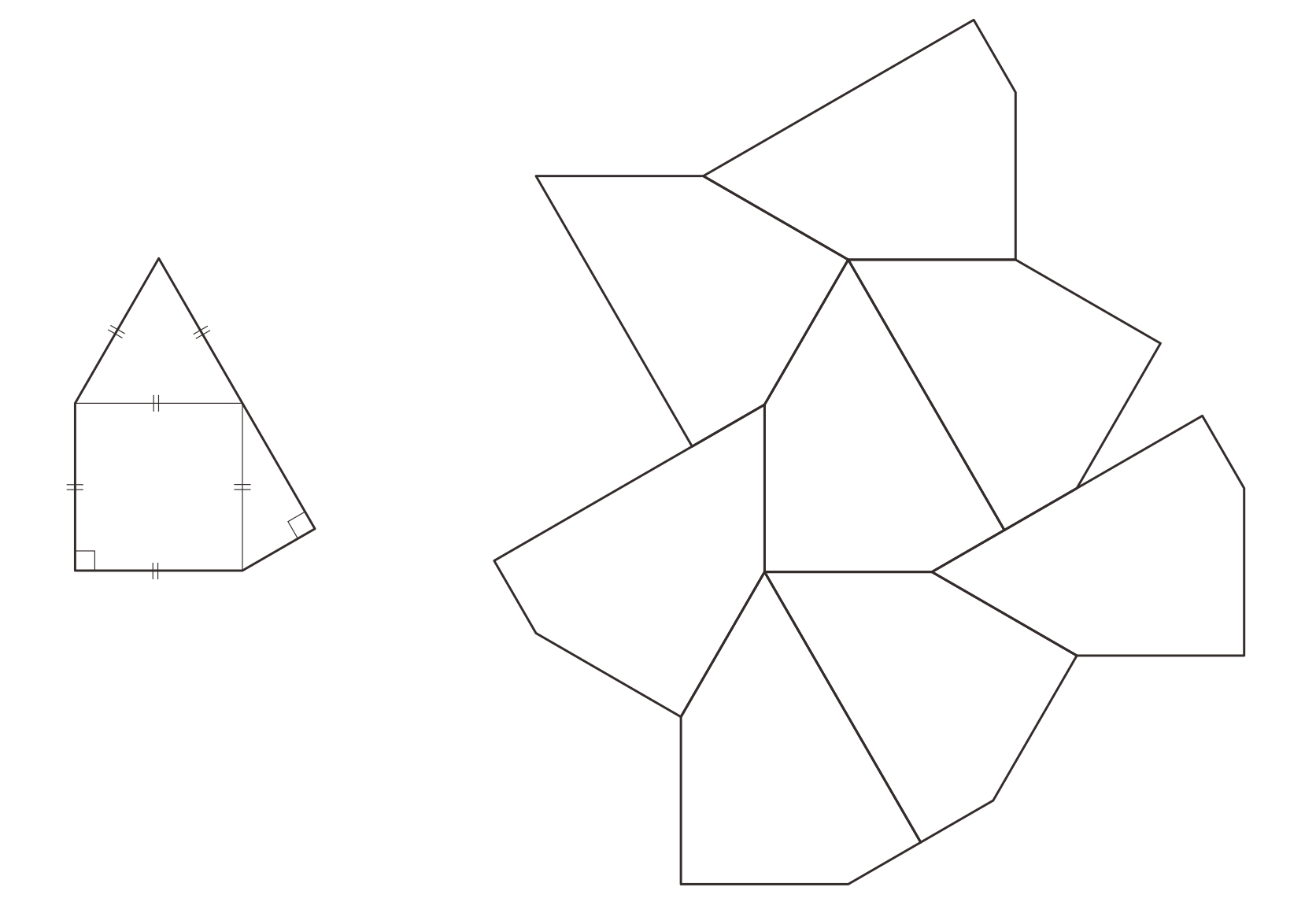} 
  \caption{{\small 
Heesch's convex pentagon with $H(T) = 1$ and its first corona. It 
seems that this convex pentagon contains a square, an equilateral 
triangle, and a $30^ \circ - 60^ \circ - 90^ \circ $ right triangle.
} 
\label{fig02}
}
\end{figure}

The Heesch number is the outcome of the presentation of a convex pentagon 
with $H(T) = 1$ (see Figure~\ref{fig02}) by Heinrich Heesch, who proposed the 
following problem in 1968~\cite{Agaoka_2005, Brass_2005, G_and_S_1987, 
Mann_2004, Wiki_HeeP}.

\bigskip
\noindent
\textbf{Heesch's Tiling Problem}. For which positive integers $k$ does there 
exist a prototile $T_{0}$ such that $T_{0}$ can be surrounded $k$ times, but not 
$k$+1 times, by tiles congruent to $T_{0}$?

\bigskip

Note that, in 1928, Lietzmann published a curvilinear tile called a spandrel 
with $H(T) = 1$~\cite{Mann_2004}.

Tiles with a finite Heesch number were known only for tiles with $H(T) = 1$ 
until 1991, but now tiles with $H(T) = 5$ are known~\cite{Mann_2004, Wiki_HeeP}. 
However, all tiles with $2 \le H(T) < \infty $ that are known so far are 
concave. The known convex polygons with a finite Heesch number are all 
polygons with $H(T) = 1$ (i.e., convex polygons that can only form the 
first corona). In 2005, Agaoka discovered the convex heptagon with 
$H(T) = 1$~\cite{Agaoka_2005}. Thereafter, Agaoka presented the problem 
``Is there a convex pentagon $T_{p}$ with $1 \le H(T) < \infty$ 
which admits an edge-to-edge corona?'' (As shown in Figure~\ref{fig02}, the first 
corona formed by the convex pentagon with $H(T) = 1$ that Heesch showed in 
1968, is non-edge-to-edge, and such a convex pentagon cannot form an 
edge-to-edge corona.) For the aforementioned Agaoka's problem, the 
solution ``$T_{p}$'' exists. Therefore, there exist convex pentagons with 
$H(T) = 1$, which admit an edge-to-edge corona. Furthermore, according to 
the same classification method of the convex pentagonal tiles (i.e., 
classification method based on the essentially different properties of 
pentagons), we find that there are infinite families of the convex pentagons 
with $H(T) = 1$. However, convex pentagons with $H(T) = 1$ should be 
classified into a finite number of families, when classified according to 
the criteria called ``Category,'' which is different from the criteria 
called ``Type'' used for convex pentagonal tiles with $H(T) = \infty $ (see 
Section~\ref{section2} for details). In this paper, we introduce the result of 
classifying convex pentagons with $H(T) = 1$, which admit edge-to-edge 
corona that we found, based on ``Category.''

\renewcommand{\figurename}{{\small Figure.}}

\renewcommand{\figurename}{{\small Figure.}}
\begin{figure}[htbp]
 \centering\includegraphics[width=13cm,clip]{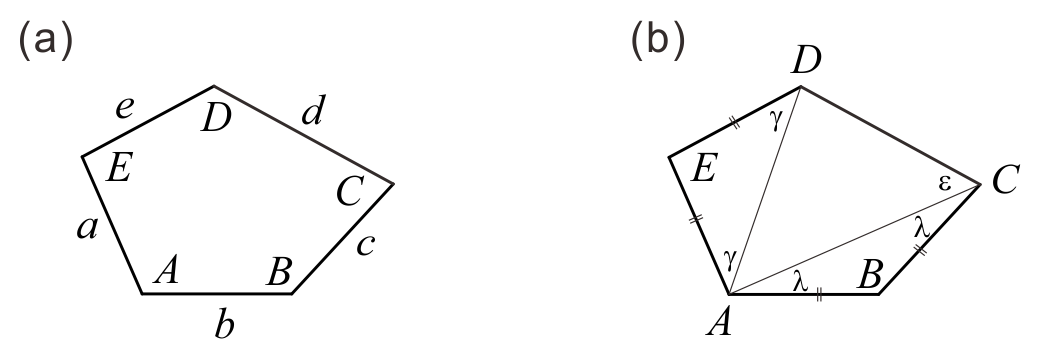} 
  \caption{{\small 
(a) Nomenclature for vertices and edges of convex pentagon. 
(b) Convex pentagon $P_{1}$ that is divided into triangles 
\textit{ABC}, \textit{ACD}, and \textit{ADE}.
} 
\label{fig03}
}
\end{figure}

\section{Convex Pentagons with $H(T) = 1$, which admit edge-to-edge corona, and 
Categories}
\label{section2}

\subsection{Preparation}
\label{subsection2.1}

Before assigning a specific Category, we need to explain the description, 
etc. Naturally, convex pentagons belonging to different Categories have 
essentially different properties, but the essential properties of convex 
pentagons within the same Category also differ. However, convex pentagons 
belonging to the same Category always have some common angle relation 
and edge relation.

In Categories, there are convex pentagons with $H(T) = 1$ and convex 
pentagons with $H(T) = \infty $. As described above, the convex pentagons 
with $H(T) = \infty $ are convex pentagons that can generate a tiling that 
infinitely covers the plane. In this paper, for clarity, when the Heesch 
number is omitted and the tile is merely described as ``convex pentagonal 
tile,'' the description refers to a convex pentagon with $H(T) = \infty $. 
For example, when we talk about ``a convex pentagonal tile belonging to Type 
6'' or ``a convex pentagonal tile of Type 6,'' we mean a convex pentagon 
belonging to (the family called) Type 6 as per the known classification 
notation, which has $H(T) = \infty $.

The convex pentagons with $H(T) = 1$ shown in this paper, possess properties 
that can form at least one pattern of edge-to-edge corona (their 
edge-to-edge corona is the first corona because $H(T) = 1)$. Some pentagons 
have either the property of forming some patterns of edge-to-edge corona or 
the property of forming non-edge-to-edge corona. In Subsection~\ref{subsection2.2}, 
we explain the properties in detail using a convex pentagon of Category 1 (details 
of such properties are omitted for convex pentagons belonging to other 
Categories).

In each Category, the ``Angle relation'' and ``Edge relation'' are the 
relations between the angles and edges that the convex pentagons in the 
Category need to satisfy, but their relations are not similar to the 
conditions of convex pentagonal tiles, as shown in Figure~\ref{fig01}. That is, the 
notations of these relations of the Category do not represent the essential 
properties of convex pentagons. The notation of the relations will be explained 
below, using a convex pentagon of Category 1. As shown in Figure~\ref{fig03}(a), 
let us label the vertices (angles) of the convex pentagon $A$, $B$, $C$, $D$, and 
$E$, as well as its edges $a$, $b$, $c$, $d$, and $e,$ in a fixed manner. As shown in 
Subsection~\ref{subsection2.2}, the convex pentagon of Category 1 has the angle 
relation ``$2A+B = 2B+E = 2D+A = 2C+A+E = 360^ \circ$,'' and the edge 
relation ``$a = b = c = e \ne d$.'' In the representation of the conditions of 
each Type of convex pentagonal tile in Figure~\ref{fig01}, for example 
when the edge condition is ``$a = b = c = d$,'' there is a possibility 
that the convex pentagon with ``$a = b = c = d = e$'' is also 
included in the Type\footnote{ For example, the convex pentagonal tiles 
of Type 7 and Type 9 in Figure~\ref{fig01} both have the edge condition of 
``$a = b = c = d$'' and have one degree of freedom except for size. As a 
result, in the family of Type 7, there exists the convex pentagonal tile 
with ``$a = b = c = d = e$.'' On the other hand, in the family of Type 9, there 
does not exist the convex pentagonal tile with ``$a = b = c = d = e$'' from 
geometrical properties~\cite{Su_and_Og_2009, Wiki_PenTP}. ``Possibility'' 
means the above-mentioned property.}. 
However, there is no such possibility with respect to the edge 
relation of convex pentagons in this section (that is, the edge relation is 
fixed). Therefore, in order to clarify, the edge relation is expressed as 
``$a = b = c = e \ne d$.'' Next, let us explain the angle relation. Here, we 
assume a convex pentagon $P_{1}$ that satisfies the relations 
``$2A+B = 2B+E = 360^ \circ$'' and ``$a = b = c = e \ne d$,'' such that 
the angle relation is looser than the convex pentagon of Category 1. 
In $P_{1}$, as shown in Figure~\ref{fig03}(b), consider two isosceles 
triangles, \textit{ABC} with base angles $\lambda $ and \textit{ADE} 
with base angles $\gamma $, and the triangle \textit{ACD}. Then, 
the interior angles of $P_{1}$ can be expressed as follows:

\begin{equation}
\label{eq1}
\left\{ {\begin{array}{l}
 A = 90^ \circ + \lambda , \\ 
 B = 180^ \circ - 2\lambda , \\ 
 C = \lambda + \varepsilon , \\ 
 D = 270^ \circ - 4\lambda - \varepsilon , \\ 
 E = 180^ \circ - 2\gamma = 4\lambda , \\ 
 \end{array}} \right.
\end{equation}

\noindent
where, the triangle \textit{ACD} and the sine theorem gives,

\[
\varepsilon = \tan ^{ - 1}\left( {\frac{\sin ^2(2\lambda )}{\cos \lambda - 
\sin (2\lambda )\cos (2\lambda )}} \right)
\]

\noindent
and we have $15.64^ \circ < \lambda < 45^ \circ $ since all interior angles 
need to be less than $180^ \circ$. Actually, the convex pentagon of Category 1 
corresponds to the case, where the convex pentagon $P_{1}$ further satisfies 
``$2D+A = 360^ \circ$.'' This implies that the convex pentagon of Category 1 
corresponds to the case, where $P_{1}$ has $\lambda = 
0.4125742...\;\mbox{rad} \approx 23.64^ \circ $, and $A \approx 113.64^ \circ ,\;
B \approx 132.72^ \circ ,\; C \approx 75.90^ \circ ,\; D \approx 123.18^ \circ , \;
E \approx 94.56^ \circ $ are the values obtained for each of the angles. The 
convex pentagon that satisfies relation (\ref{eq1}) and ``$2D+A = 360^ \circ$'' 
also satisfies ``$2C+A+E = 360^ \circ$.'' Therefore, with respect to the angle 
relation (angle condition) of the convex pentagon of Category 1, ``$2A+B = 
2B+E = 2D+A = 360^ \circ$'' is sufficient, but in this paper, we describe the angle 
relation as ``$2A+B = 2B+E = 2D+A = 2C+A+E = 360^ \circ$.'' This notation is a 
collection of the vertices' concentrating methods (relations where the sum of interior 
angles is $360^ \circ$) used for forming corona exemplified.

As shown in the example of the convex pentagon of Category 1, the value of 
each interior angle of the convex pentagon of each Category is described in 
the table. The exact values (obtained from mathematical formulas) of each 
interior angle can be derived from the relations between angles and edges as 
in the above example, but it is difficult to describe all of them. Hence, 
for most interior angle values in the tables in each Category, only the 
approximate values are indicated.

The item ``Heesch number'' shown in each Category is the value of $H(T)$ of the 
convex pentagons included in the Category.

 Let us call the multi-set of vertices of pentagons a \textit{spot}. If the sum of the 
interior angles at the vertices in the multi-set is equal to $360^ \circ$, then 
the spot that can be concentrated by edge-to-edge contact is called 
\textit{EEC-spot}. For example, as show in Figure~\ref{fig05}, the spot ``$2C+A+E = 360^ \circ$'' 
of Category 1 can be concentrated by non-edge-to-edge contact, but it can also be 
concentrated by edge-to-edge contact. Therefore, the spot ``$2C+A+E = 360^ \circ$'' 
is the EEC-spot. On other hand, the spot that cannot be concentrated by 
edge-to-edge contact is called \textit{NEEC-spot}.

\subsection{Classification result by Category}
\label{subsection2.2}

In this subsection, seventeen families of Categories and the convex 
pentagons belonging to each Category are introduced.

\bigskip\bigskip
\noindent
\textbf{Category 1}

\begin{figure}[htbp]
 \centering\includegraphics[width=14.5cm,clip]{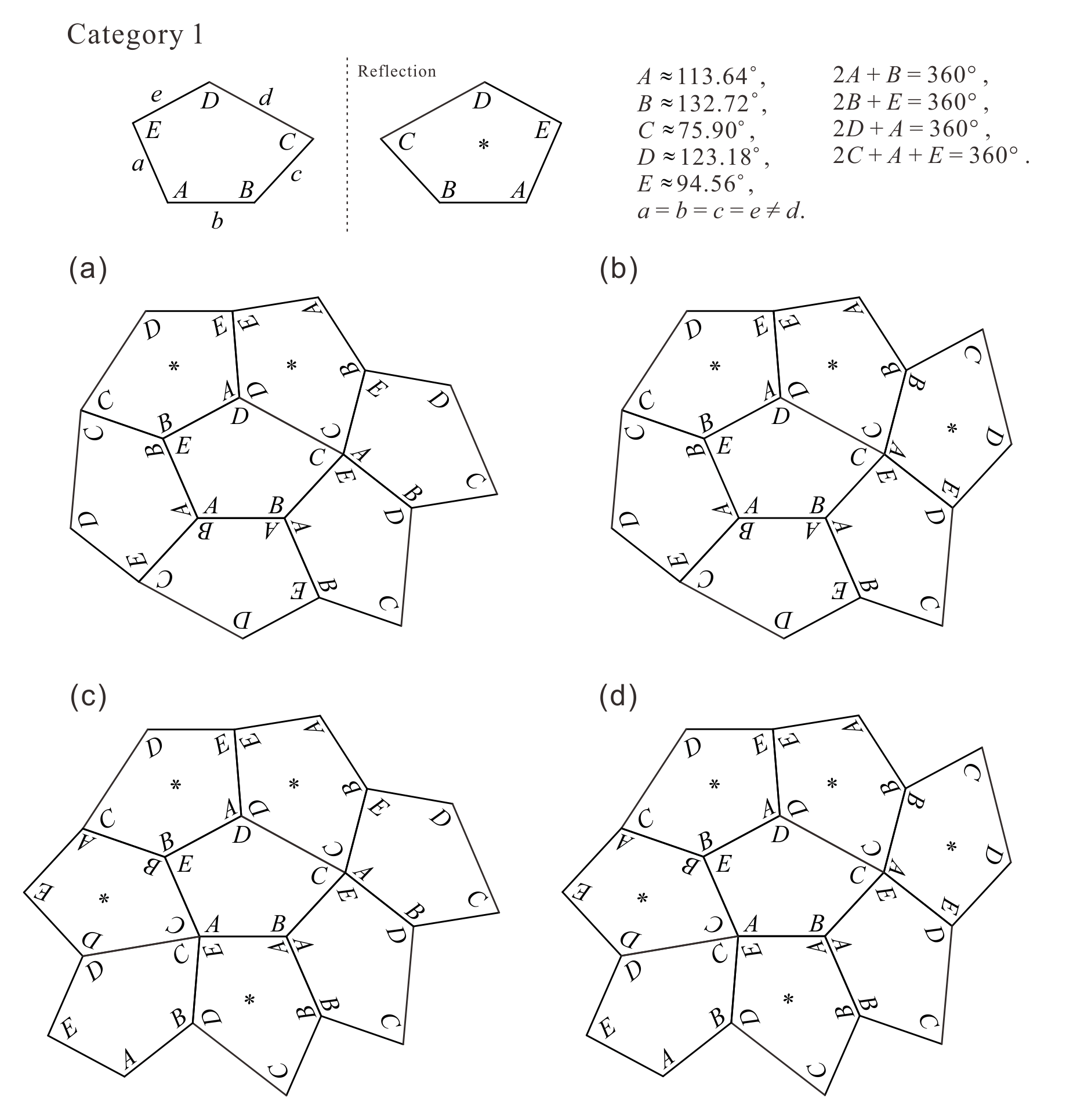} 
  \caption{{\small 
Convex pentagon of Category 1 and edge-to-edge coronas by the 
pentagons.
} 
\label{fig04}
}
\end{figure}

\begin{figure}[htbp]
 \centering\includegraphics[width=14.5cm,clip]{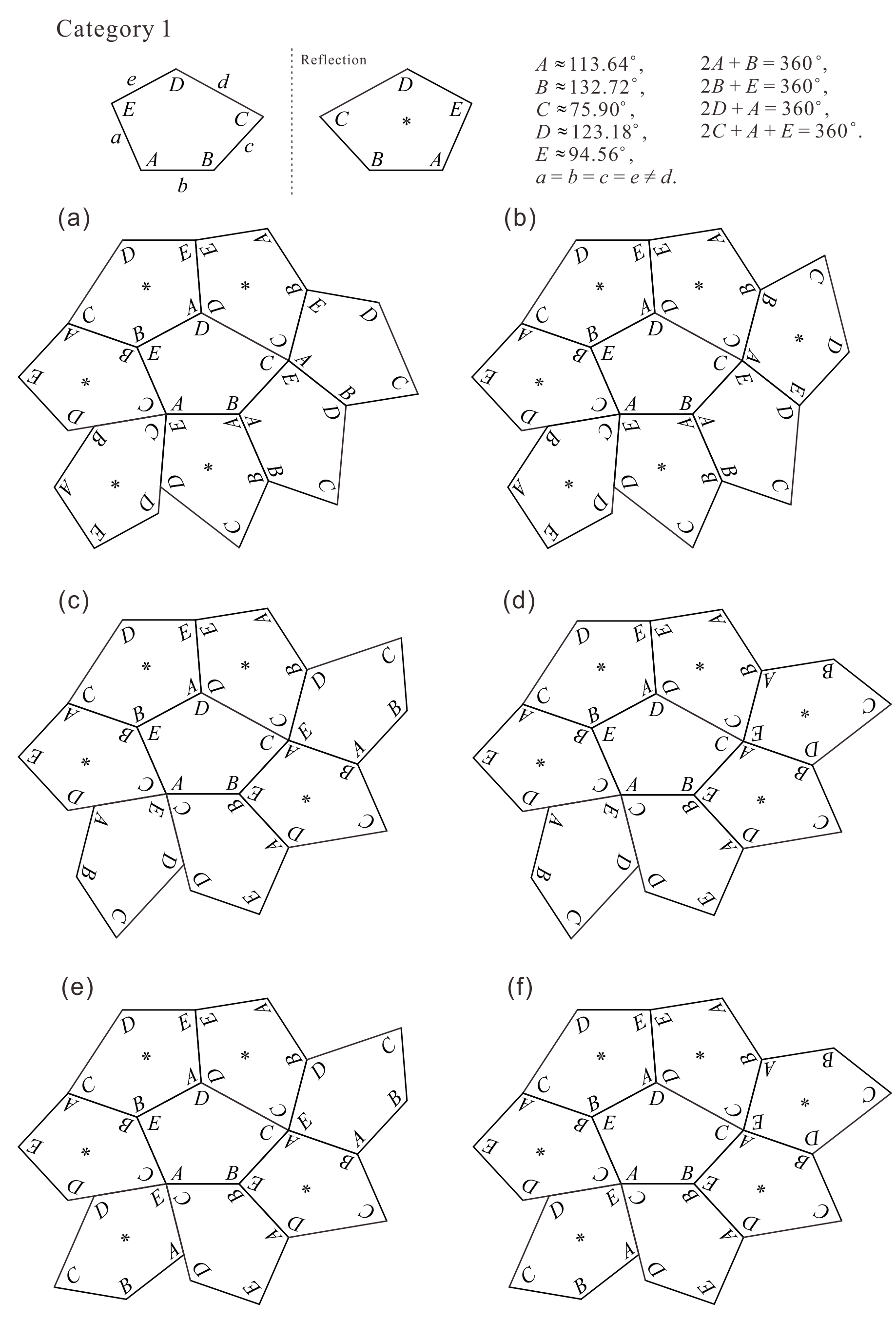} 
  \caption{{\small 
Convex pentagon of Category 1 and non-edge-to-edge coronas by the 
pentagons.
} 
\label{fig05}
}
\end{figure}

\begin{description}
 \setlength{\itemindent}{-10pt}
 \setlength{\itemsep}{-3pt} 
\item[Angle relation:] $2A+B = 2B+E = 2D+A = 2C+A+E = 360^ \circ$.

\item[Edge relation:] $a = b = c = e \ne d$.

\item[Heesch number:] $H(T) = 1$.

\item[Corresponding Table and Figures:] Table~\ref{tab01}, and Figures~\ref{fig04} and \ref{fig05}.
\end{description}

\begin{table}[!h]
 \begin{center}
{\small
\caption[Table 1]{Value and arrangement of vertices of convex pentagon of Category 1}
\label{tab01}
}

\
{\footnotesize
\begin{tabular}
{rrrrr|rrrrr}
\hline
\multicolumn{5}{c|}{\raisebox{-1.75ex}[0.5cm][0.5cm]
{\footnotesize \shortstack{ Value of interior angle \\(degree) }}  } & 
\multicolumn{5}{c}{\raisebox{-1.75ex}[0.5cm][0.5cm]
{\footnotesize \shortstack{ Example of arrangement around each \\vertex (counterclockwise) }} }  \\

$A$& 
$B$& 
$C$& 
$D$& 
$E$& 
$A$& 
$B$& 
$C$& 
$D$& 
$E$ \\
\hline

113.64 & 
132.72 & 
75.90 & 
123.18 & 
94.56 & 
\textit{AAB}& 
\textit{BAA}& 
\textit{CEAC}& 
\textit{DDA}& 
\textit{EBB} \\
\hline

\end{tabular}
}
\end{center}
\end{table}

\noindent
\textbf{Remarks. }
The shapes of the edge-to-edge corona that can be formed 
by a convex pentagon of Category 1 are the four patterns shown in Figure~\ref{fig04}. 
For example, consider ``\textit{CEAC}'' as the counterclockwise arrangement around 
vertex $C$ as shown in Table~\ref{tab01}. It means that the vertices concentrate at 
a point in the order ``$C \to E \to A \to C$'' in the counterclockwise direction on 
the vertex $C$ of the centrally convex pentagon as shown in Figure~\ref{fig04}. On 
the other hand, the convex pentagons of Category 1 can also form non-edge-to-edge 
coronas, which can be seen in the six patterns in Figure~\ref{fig05}.

\begin{figure}[htbp]
 \centering\includegraphics[width=14.5cm,clip]{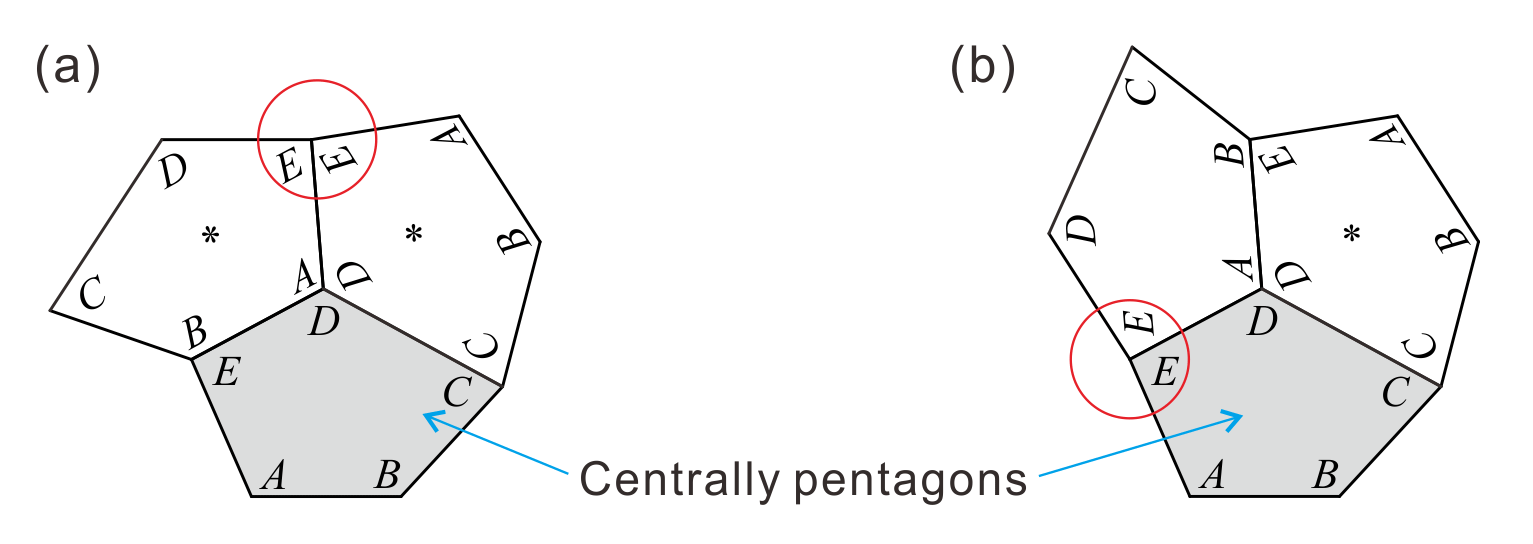} 
  \caption{{\small 
Patterns where vertex $D$ of convex pentagon of Category 1 is 
$2D+A=360^ \circ$.
} 
\label{fig06}
}
\end{figure}

\begin{figure}[htbp]
 \centering\includegraphics[width=14.5cm,clip]{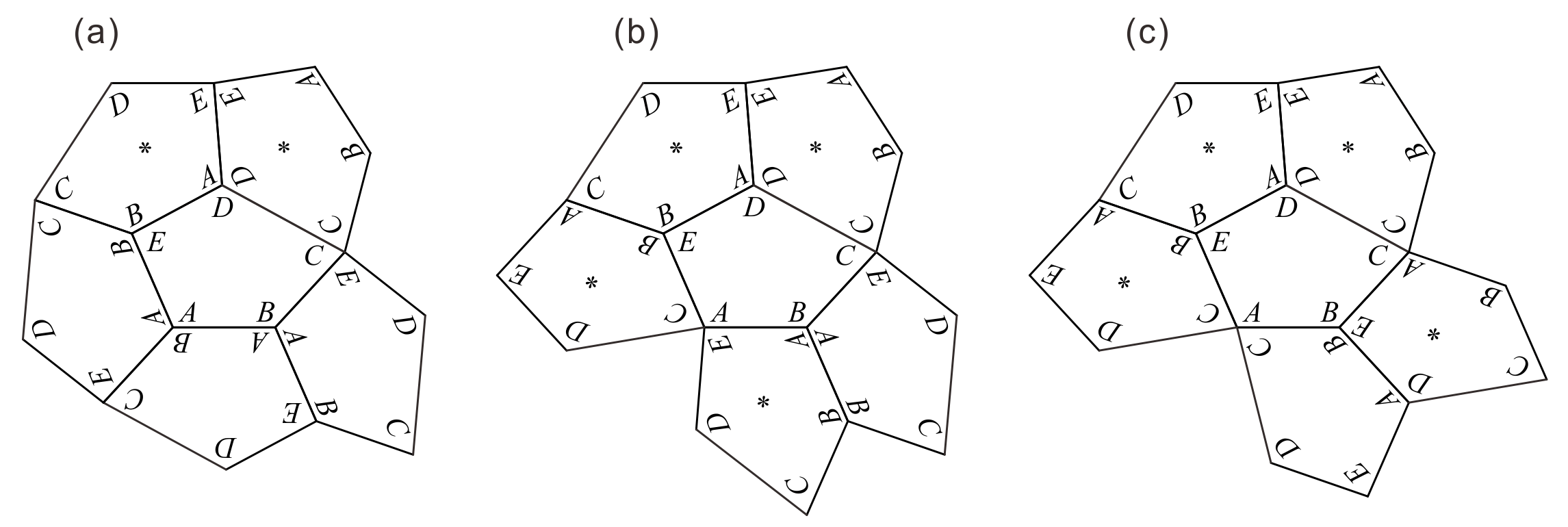} 
  \caption{{\small 
Based patterns when the convex pentagon of Category 1 forms the 
first corona.
} 
\label{fig07}
}
\end{figure}

Next, we show that the convex pentagon of Category 1 cannot have 
$H(T) \geq 2$ (this explanation is omitted in other Categories). In order 
for a convex pentagon to generate tiling, it is necessary that two or more of its 
vertices are concentrated, and the sum of internal angles in the concentration 
of vertices must be $180^ \circ$ or $360^ \circ$. ``$2A+B = 2B+E = 2D+A = 
2C+A+E = 360^ \circ$'' is the relation that the sum of the interior angles of the 
convex pentagon of Category 1 is $180^ \circ$ or $360^ \circ$. Therefore, 
vertex $D$ necessarily uses the relation ``$2D+A = 360^ \circ$,'' and the 
centrally convex pentagon makes edge-to-edge contact with the other convex 
pentagons, since there is no relation that the sum of the interior angles is 
$180^ \circ$. There are only two patterns (see Figure~\ref{fig06}) where the 
vertex $D$ of the centrally convex pentagon is $2D+A = 360^ \circ$. 
In Figure~\ref{fig06}, the places where the two vertices $E$ are surrounded 
by the red circle, cannot make a relation of $360^ \circ$, and the one that 
can make the first corona out of the two patterns in Figure~\ref{fig06}, 
is the pattern of Figure~\ref{fig06}(a). Based on the pattern of Figure~\ref{fig06}(a), 
three patterns in Figure~\ref{fig07} can be created. The patterns of 
Figures~\ref{fig04}(a) and \ref{fig04}(b) are created from the patterns of 
Figure~\ref{fig07}(a), the patterns of Figures~\ref{fig04}(c), \ref{fig04}(d), 
\ref{fig05}(a), and \ref{fig05}(b) are created from the patterns of 
Figure~\ref{fig07}(b), and the patterns of Figures~\ref{fig05}(c), \ref{fig05}(d), 
\ref{fig05}(e), and \ref{fig05}(f) are created from the patterns 
of Figure~\ref{fig07}(c). Since in each pattern, there is a place where two vertices 
$E$ are concentrated at the boundaries of the first corona of Figures~\ref{fig04} 
and \ref{fig05}, $H(T) \geq 2 $ is not permitted for the convex pentagon of Category 1.

\renewcommand{\figurename}{{\small Figure.}}

\renewcommand{\figurename}{{\small Figure.}}

\renewcommand{\figurename}{{\small Figure.}}

\renewcommand{\figurename}{{\small Figure.}}

\bigskip\bigskip
\noindent
\textbf{Category 2}

\begin{figure}[htbp]
 \centering\includegraphics[width=14.5cm,clip]{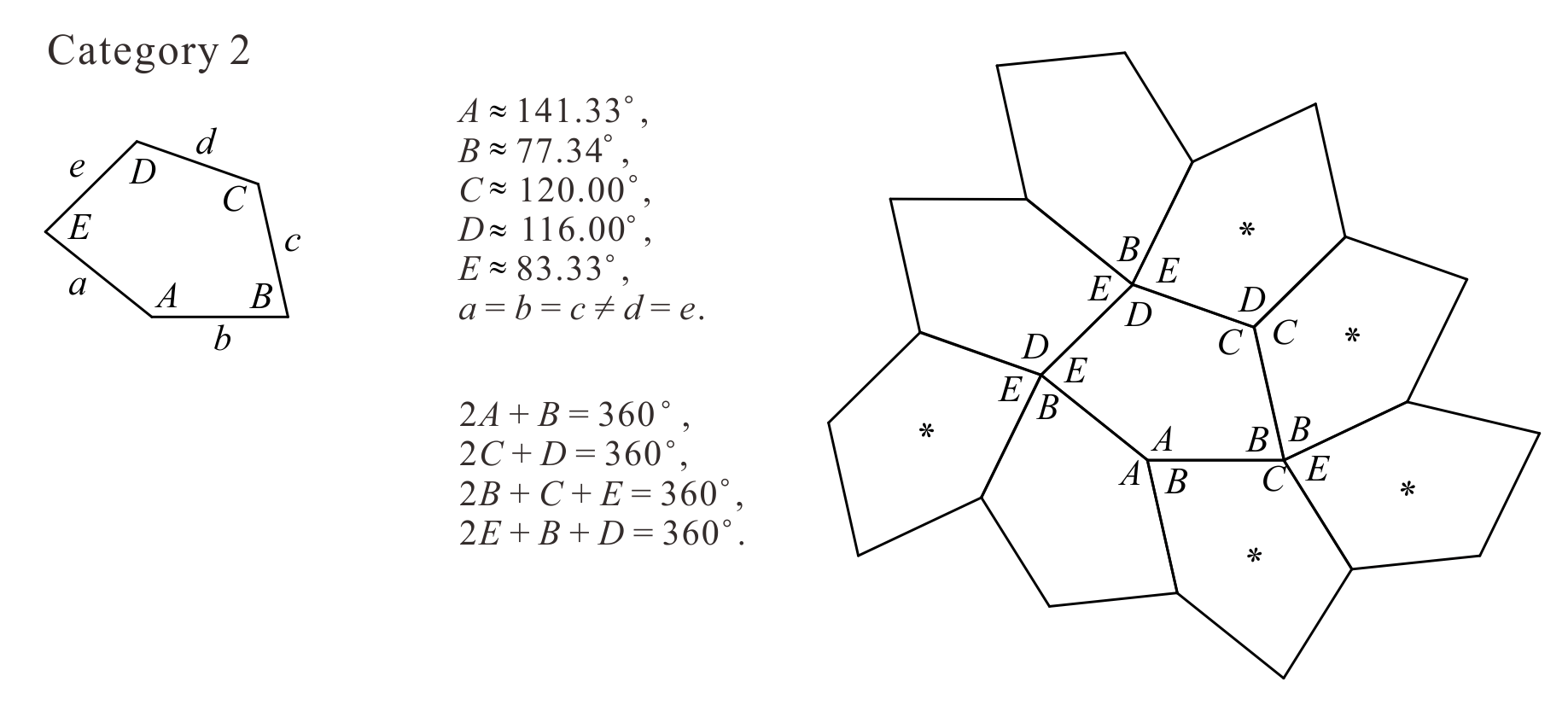} 
  \caption{{\small 
Convex pentagon of Category 2 and example of edge-to-edge corona 
by the pentagons.
} 
\label{fig08}
}
\end{figure}

\begin{description}
 \setlength{\itemindent}{-10pt}
 \setlength{\itemsep}{-3pt} 
\item[Angle relation:] $2A+B = 2C+D = 2B+C+E = 2E+B+D = 360^ \circ$.

\item[Edge relation:] $a = b = c \ne d = e$.

\item[Heesch number:] $H(T) = 1$.

\item[Corresponding Table and Figure:] Table~\ref{tab02} and Figure~\ref{fig08}.
\end{description}

\begin{table}[!h]
 \begin{center}
{\small
\caption[Table 2]{Value and arrangement of vertices of convex pentagon of Category 2}
\label{tab02}
}
\
{\footnotesize
\begin{tabular}
{rrrrr|rrrrr}
\hline
\multicolumn{5}{c|}{\raisebox{-1.75ex}[0.5cm][0.5cm]
{\small \shortstack{ Value of interior angle \\(degree) }}  } & 
\multicolumn{5}{c}{\raisebox{-1.75ex}[0.5cm][0.5cm]
{\small \shortstack{ Example of arrangement around each \\vertex (counterclockwise) }} }  \\

$A$& 
$B$& 
$C$& 
$D$& 
$E$& 
$A$& 
$B$& 
$C$& 
$D$& 
$E$ \\
\hline

141.33 & 
77.34 & 
122.00 & 
116.00 & 
83.33 & 
\textit{AAB}& 
\textit{BCEB}& 
\textit{CCD}& 
\textit{DEBE}& 
\textit{EDEB} \\
\hline

\end{tabular}
}
\end{center}
\end{table}

\renewcommand{\figurename}{{\small Figure.}}

\bigskip\bigskip
\noindent
\textbf{Category 3}

\begin{figure}[htbp]
 \centering\includegraphics[width=14.5cm,clip]{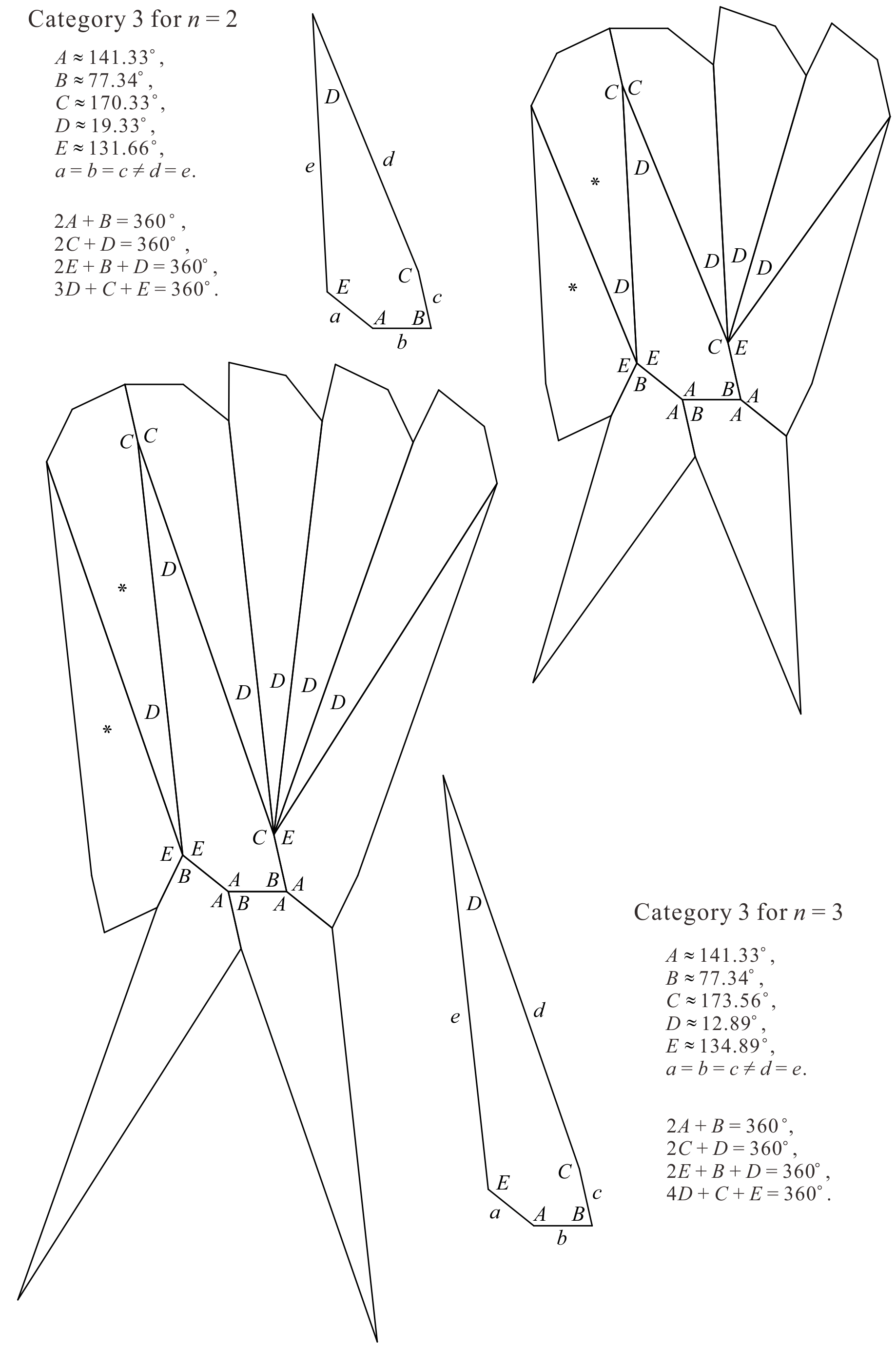} 
  \caption{{\small 
Convex pentagons of Category 3 and examples of edge-to-edge coronas 
by the pentagons.
} 
\label{fig09}
}
\end{figure}

\begin{description}
 \setlength{\itemindent}{-10pt}
 \setlength{\itemsep}{-3pt} 
\item[Angle relation:] $2A+B = 2C+D = 2E+B+D = (n+1)\times D+C+E = 360^ \circ$ where $n$ is an 
integer of one or more.

\item[Edge relation:] $a = b = c \ne d = e$.

\item[Heesch number:]  $H(T) = \infty$ for $n=1$ or $H(T) = 1$ for $n > 1$.

\item[Corresponding Table and Figure:] Table~\ref{tab03} and Figure~\ref{fig09}.
\end{description}

\begin{table}[!h]
 \begin{center}
{\small
\caption[Table 3]{Value and arrangement of vertices of convex pentagon of Category 3}
\label{tab03}
}
\
{\footnotesize
\begin{tabular}
{c|c|rrrrr|rrrrr}
\hline
\raisebox{-1.50ex}[0cm][0cm]{$n$}& 
\raisebox{-1.50ex}[0cm][0cm]{$H(T)$}& 
\multicolumn{5}{c|}{\raisebox{-1.75ex}[0.5cm][0.5cm]
{\footnotesize \shortstack{ Value of interior angle \\(degree) }}  } & 
\multicolumn{5}{c}{\raisebox{-1.75ex}[0.5cm][0.5cm]
{\footnotesize \shortstack{ Example of arrangement around each \\vertex (counterclockwise) }} }  \\
 & 
 & 

$A$& 
$B$& 
$C$& 
$D$& 
$E$& 
$A$& 
$B$& 
$C$& 
$D$& 
$E$ \\
\hline
 
1& 
$\infty $& 
141.33 & 
77.34& 
160.67& 
38.67& 
122.00& 
\multicolumn{5}{c}{Convex Pentagonal tile belonging to Type 6}  \\

\hline

2& 
1& 
141.33 & 
77.34 & 
170.33 & 
19.33 & 
131.66 & 
\textit{AAB}& 
\textit{BAA}& 
\textit{CEDDD}& 
\textit{DCC}& 
\textit{EDEB} \\
\hline

3& 
1& 
141.33 & 
77.34 & 
173.56 & 
12.89 & 
134.89 & 
\textit{AAB}& 
\textit{BAA}& 
\textit{CEDDDD}& 
\textit{DCC}& 
\textit{EDEB} \\
\hline

4& 
1& 
141.33 & 
77.34 & 
175.17 & 
9.67 & 
136.50 & 
\textit{AAB}& 
\textit{BAA}& 
\textit{CEDDDDD}& 
\textit{DCC} & 
\textit{EDEB} \\
\hline

...& 
1& 
...& 
...& 
...& 
...& 
...& 
... & 
... & 
... & 
... & 
...  \\

\hline
\end{tabular}
}
\end{center}
\end{table}

\noindent
\textbf{Remarks. }Firstly, a convex pentagon satisfying ``$2A + B = 2C + D = 
D + C + E = 360^ \circ , a = b = c \ne d = e$'' is geometrically impossible 
(i.e., the case of $n=0$ does not exist). 
Secondly, the convex pentagon for $n = 1$ is the convex pentagonal tile 
belonging to Type 6 (in this case, EEC-spots ``$ B+C+E = 3D+2E = 360^ \circ$'' and 
NEEC-spot ``$2A+2D = 360^ \circ$'' also hold). Therefore, the convex pentagon of 
Category 3 is the convex pentagon with $H(T) = \infty $ if and only if $n = 1$. 
Then, for $n \ge 2$, EEC-spot ``$(2n+1) \times D+2E = 360^ \circ$'' and 
NEEC-spot ``$2n \times D+2A = 360^ \circ$'' should also hold.

\begin{figure}[htbp]
 \centering\includegraphics[width=13cm,clip]{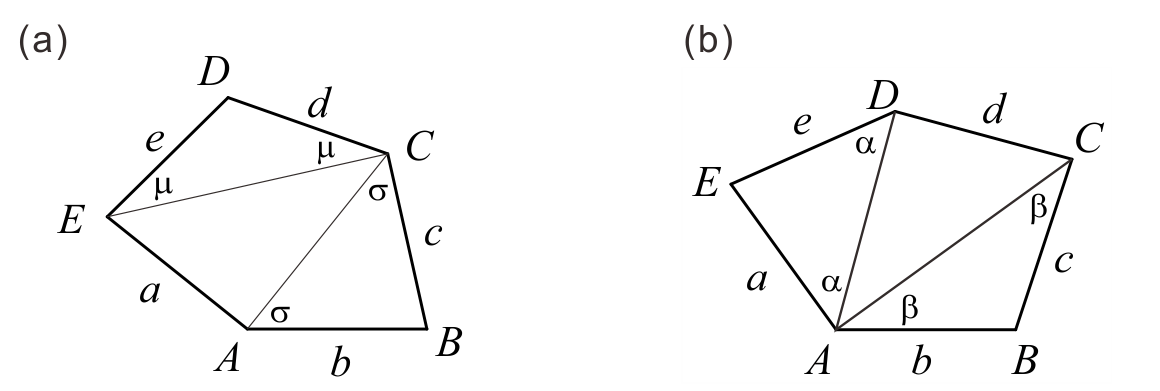} 
  \caption{{\small 
Convex pentagons $P_{2}$ and $P_{3.}$ (a)$ P_{2}$ that satisfies 
``$2A+B = 2C+D = 360^ \circ, a = b = c \ne d = e$.'' (b) $P_{3}$ that satisfies 
``$2A+B = 2D+E = 360^ \circ, a = b = c = e \ne d$.''
} 
\label{fig10}
}
\end{figure}

Let $P_{2}$ be a convex pentagon that satisfies ``$2A+B = 2C+D = 360^ \circ, 
a = b = c \ne d = e$.'' In $P_{2}$, as shown in Figure~\ref{fig10}(a), consider two 
isosceles triangles, \textit{ECD} with base angles $\mu $ and \textit{ABC} with base 
angles $\sigma $, then the interior angles of $P_{2}$ can be expressed as follows:

\begin{equation}
\label{eq2}
\left\{ {\begin{array}{l}
 A = 90^ \circ + \sigma , \\ 
 B = 180^ \circ - 2\sigma , \\ 
 C = 90^ \circ + \mu , \\ 
 D = 180^ \circ - 2\mu , \\ 
 E = \mu + \sigma , \\ 
 \end{array}} \right.
\end{equation}

\noindent
where $0^ \circ < \mu < 90^ \circ $ from $D < 180^ \circ $ and $C < 180^ 
\circ $, and

\begin{equation}
\label{eq3}
\sigma = \sin ^{ - 1}\left( {\frac{ - 1 + \sqrt {17} }{4}} \right) \approx 
0.895907\;\mbox{rad} \approx 51.3317^ \circ 
\end{equation}

\noindent
using the triangle \textit{ACE} and the sine theorem.

If $P_{2}$ has the relation ``$(n + 1)\times D + C + E = 360^ \circ $,'' from 
$(n + 1)\times D + C + E = 360^ \circ $ and (\ref{eq2}), we get

\begin{equation}
\label{eq4}
n = \frac{270^ \circ - \sigma - 2\mu }{180^ \circ - 2\mu } - 1.
\end{equation}

\noindent
From $270^ \circ - \sigma \approx 218.66^ \circ $, 
$270^ \circ - \sigma - 2\mu > 180^ \circ - 2\mu > 0^ \circ $ holds 
for $0^ \circ < \mu < 90^ \circ $. Therefore, $n > 0$ for 
$0^ \circ < \mu < 90^ \circ $. From (\ref{eq2}), (\ref{eq3}) and (\ref{eq4}), 
as $\mu \to 90^ \circ $, we see that $n \to \infty $. Thus, there are 
infinite convex pentagons with $H(T) = 1$ belonging to Category 3.

\renewcommand{\figurename}{{\small Figure.}}

\bigskip\bigskip
\bigskip\bigskip
\bigskip\bigskip
\bigskip\bigskip
\bigskip\bigskip

\bigskip\bigskip
\noindent
\textbf{Category 4}

\begin{figure}[htbp]
 \centering\includegraphics[width=14.5cm,clip]{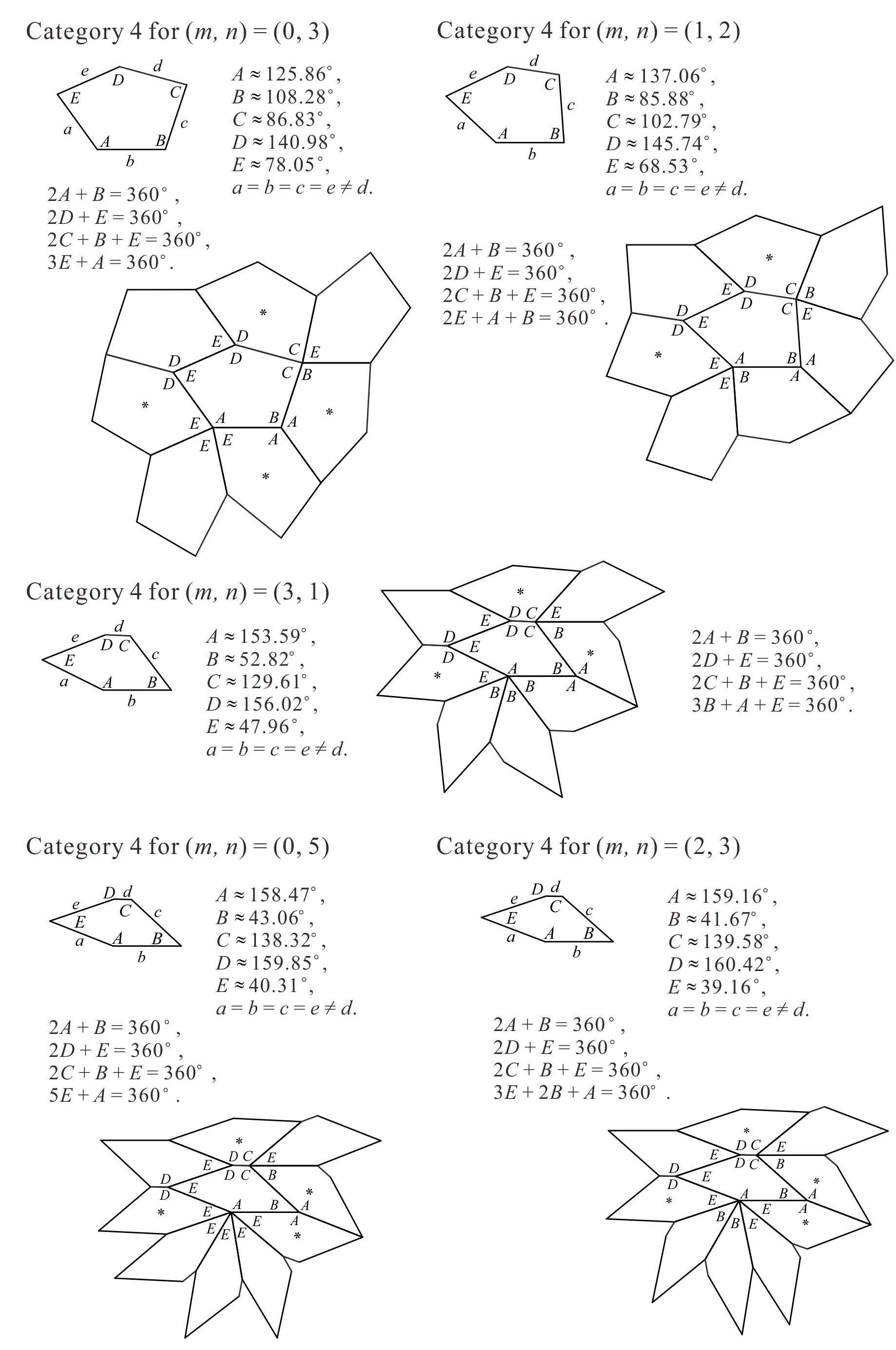} 
  \caption{{\small 
Convex pentagons of Category 4 and examples of edge-to-edge coronas 
by the pentagons.
} 
\label{fig11}
}
\end{figure}

\begin{description}
 \setlength{\itemindent}{-10pt}
 \setlength{\itemsep}{-3pt} 
\item[Angle relation:] $2A+B = 2D+E = 2C+B+E= m \times B+n \times E+A = 360^ \circ$ 
where $m$ is an integer of zero or more, $n$ is an integer of one or more,   
and $m + n \ge 3$.

\item[Edge relation:] $a = b = c = e \ne d$.

\item[Heesch number:] $H(T) = 1$ for $(m,n) \ne (2,\;1)$ or $H(T) = \infty $ for 
$(m,n) = (2,\;1)$.

\item[Corresponding Table and Figure:] Table~\ref{tab04} and Figure~\ref{fig11}.
\end{description}

\begin{table}[!h]
 \begin{center}
{\small
\caption[Table 4]{Value and arrangement of vertices of convex pentagon of Category 4}
\label{tab04}
}
\
{\footnotesize
\begin{tabular}
{c|c|rrrrr|rrrrr}
\hline
\raisebox{-1.50ex}[0cm][0cm]{$m, n$}& 
\raisebox{-1.50ex}[0cm][0cm]{$H(T)$}& 
\multicolumn{5}{c|}{\raisebox{-1.75ex}[0.5cm][0.5cm]
{\footnotesize \shortstack{ Value of interior angle \\(degree) }}  } & 
\multicolumn{5}{c}{\raisebox{-1.75ex}[0.5cm][0.5cm]
{\footnotesize \shortstack{ Example of arrangement around each \\vertex (counterclockwise) }} }  \\
 & 
 &
 
$A$& 
$B$& 
$C$& 
$D$& 
$E$& 
$A$& 
$B$& 
$C$& 
$D$& 
$E$ \\
\hline

0, 3& 
1& 
125.86 & 
108.28 & 
86.83 & 
140.98 & 
78.05 & 
\textit{AEEE}& 
\textit{BAA}& 
\textit{CBEC}& 
\textit{DDE}& 
\textit{EDD} \\
\hline

1, 2& 
1& 
137.06 & 
85.88 & 
102.79 & 
145.74 & 
68.53 & 
\textit{AEEB}& 
\textit{BAA}& 
\textit{CEBC}& 
\textit{DDE}& 
\textit{EDD} \\
\hline

2, 1& 
$\infty $& 
141.33 & 
77.34 & 
109.33 & 
148.00 & 
64.00 & 
\multicolumn{5}{c}{Convex Pentagonal tile belonging to Type 9}  \\
\hline

0, 4& 
1& 
150.55 & 
58.90 & 
124.37 & 
153.82 & 
52.36 & 
\textit{AEEEE}& 
\textit{BAA}& 
\textit{CBEC}& 
\textit{DDE}& 
\textit{EDD} \\
\hline

1, 3& 
1& 
151.79 & 
56.42 & 
126.49 & 
154.70 & 
50.60 & 
\textit{AEEBE}& 
\textit{BAA}& 
\textit{CEBC}& 
\textit{DDE}& 
\textit{EDD} \\
\hline

2, 2& 
1& 
152.78 & 
54.45 & 
128.19 & 
155.42 & 
49.16 & 
\textit{AEBBE}& 
\textit{BAA}& 
\textit{CBEC}& 
\textit{DDE}& 
\textit{EDD} \\
\hline

3, 1& 
1& 
153.59 & 
52.82 & 
129.61 & 
156.02 & 
47.96 & 
\textit{AEBBB}& 
\textit{BAA}& 
\textit{CBEC}& 
\textit{DDE}& 
\textit{EDD} \\
\hline

0, 5& 
1& 
158.47 & 
43.06 & 
138.32 & 
159.85 & 
40.31 & 
\textit{AEEEEE}& 
\textit{BAA}& 
\textit{CBEC}& 
\textit{DDE}& 
\textit{EDD} \\
\hline

1, 4& 
1& 
158.83 & 
42.33 & 
138.98 & 
160.15 & 
39.71 & 
\textit{AEBEEE}& 
\textit{BAA}& 
\textit{CBEC}& 
\textit{DDE}& 
\textit{EDD} \\
\hline

2, 3& 
1& 
159.16 & 
41.67 & 
139.58 & 
160.42 & 
39.16 & 
\textit{AEBBEE}& 
\textit{BAA}& 
\textit{CBEC}& 
\textit{DDE}& 
\textit{EDD} \\
\hline

3, 2& 
1& 
159.46 & 
41.07 & 
140.13 & 
160.67 & 
38.66 & 
\textit{AEBBBE}& 
\textit{BAA}& 
\textit{CBEC}& 
\textit{DDE}& 
\textit{EDD} \\
\hline

4, 1& 
1& 
159.74 & 
40.52 & 
140.64 & 
160.90 & 
38.20 & 
\textit{AEBBBB}& 
\textit{BAA}& 
\textit{CBEC}& 
\textit{DDE}& 
\textit{EDD} \\
\hline

...& 
1& 
...& 
...& 
...& 
...& 
...& 
... & 
... & 
... & 
... & 
...  \\
\hline

\end{tabular}
}
\end{center}
\end{table}

\noindent
\textbf{Remarks}. As for the convex pentagon satisfying ``$2A + B = 2D + E = 
360^ \circ , a = b = c = e \ne d$,'' if the convex pentagon also has the 
relation ``$2B+E+A = 360^ \circ$,'' then it is the convex pentagonal tile belonging 
to Type 9 (in this case, EEC-spot ``$3B+2E = 360^ \circ$'' also holds). Therefore, 
the convex pentagon of Category 4 is the convex pentagon with $H(T) = \infty $ 
if and only if $(m, n) = (2,\;1)$. Note that, in the cases of $m \ne 0$, EEC-spot 
``$ (2m - 1)\times B + 2n\times E = 360^ \circ $'' also holds.

Let us assume a convex pentagon $P_{3}$ that satisfies ``$2A + B = 2D + E = 360^ \circ , 
a = b = c = e \ne d$.'' As shown in Figure~\ref{fig10}(b), $P_{3}$ contains two isosceles 
triangles, \textit{ADE} with base angles $\alpha $ and \textit{ABC} with base 
angles $\beta $. Therefore, we see that the interior angles of $P_{3}$ can be 
expressed as follows:

\begin{equation}
\label{eq5}
\left\{ {\begin{array}{l}
 A = 90^ \circ + \beta , \\ 
 B = 180^ \circ - 2\beta , \\ 
 C = \alpha + \beta , \\ 
 D = 90^ \circ + \alpha , \\ 
 E = 180^ \circ - 2\alpha , \\ 
 \end{array}} \right.
\end{equation}

\noindent
where, we can see from the triangle \textit{ACD} and the sine theorem,

\begin{equation}
\label{eq6}
\beta = \cos ^{ - 1}\left( {\frac{\cos \alpha }{\sin \alpha }} \right),
\end{equation}

\noindent
and $45^ \circ < \alpha < 90^ \circ $ ($0^ \circ < \beta < 90^ \circ $) 
since $B < 180^ \circ $ and $\beta > 0^ \circ $.

If $P_{3}$ has the relation ``$m\times B + n\times E + A = 360^ \circ $,'' 
then the integers $m$ and $n$ have relations ``$m + n \ge 3,\; m \ge 0,\;n \ge 1$'' 
from the geometric nature. Therefore, when $m = 0$, from $n\times E + A = 360^ \circ $ 
and (\ref{eq5}), we have

\begin{equation}
\label{eq7}
n = \frac{270^ \circ - \beta }{180^ \circ - 2\alpha }.
\end{equation}

\noindent
Then, $270^ \circ - \beta > 180^ \circ - \alpha > 0^ \circ $ for 
$45^ \circ < \alpha < 90^ \circ $, and $\beta = 0^ \circ $ for $\alpha = 90^ \circ $. 
Therefore, from (\ref{eq5}), (\ref{eq6}) and (\ref{eq7}), as $\alpha \to 90^ \circ $, 
we see that $n \to \infty $. Thus, there are infinite convex pentagons with $H(T) = 1$ 
belonging to Category 4.

\renewcommand{\figurename}{{\small Figure.}}

\renewcommand{\figurename}{{\small Figure.}}

\bigskip\bigskip
\noindent
\textbf{Category 5}

\begin{figure}[htbp]
 \centering\includegraphics[width=14.5cm,clip]{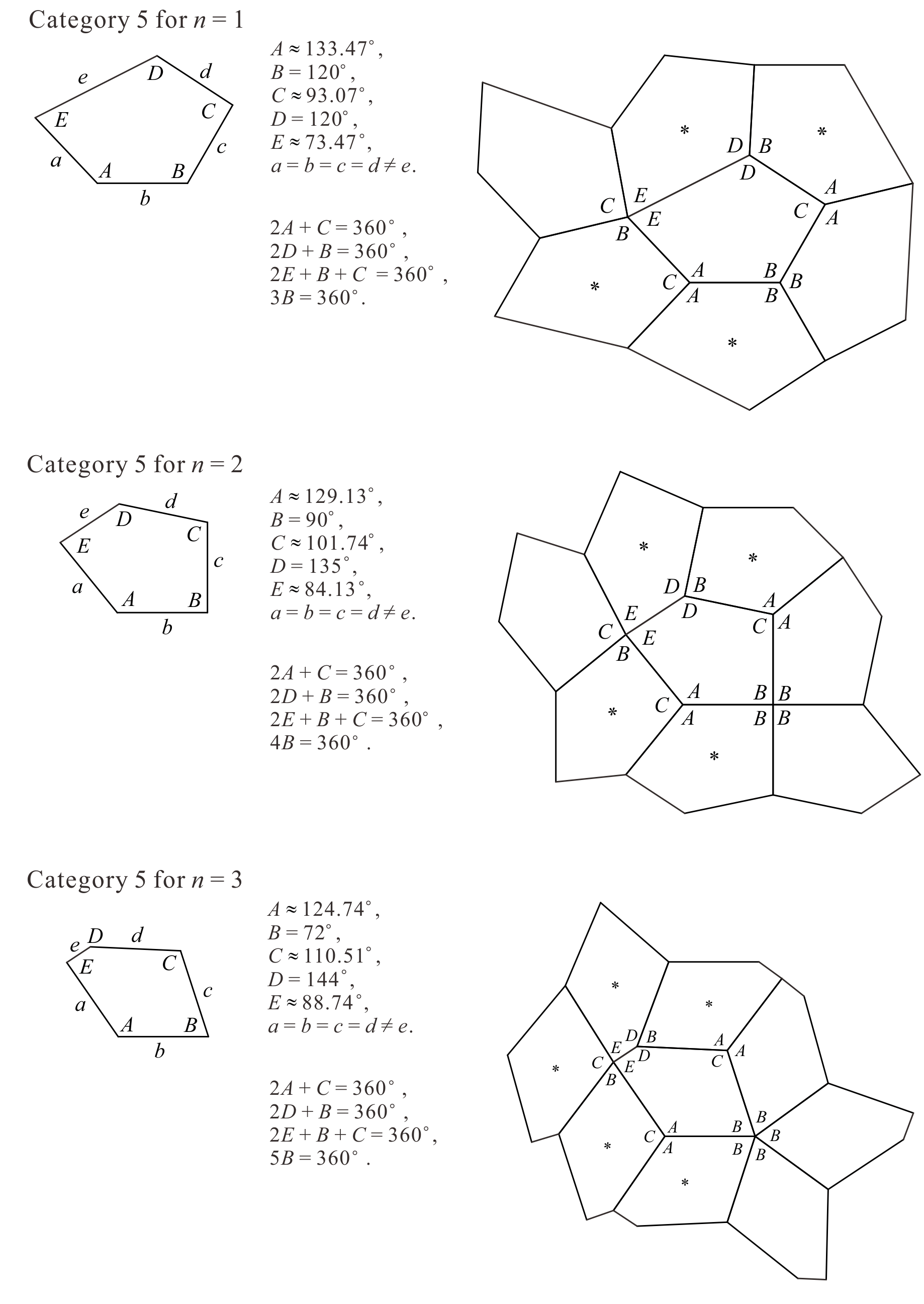} 
  \caption{{\small 
Convex pentagons of Category 5 and examples of edge-to-edge coronas 
by the pentagons.
} 
\label{fig12}
}
\end{figure}

\begin{description}
 \setlength{\itemindent}{-10pt}
 \setlength{\itemsep}{-3pt} 
\item[Angle relation:] $2A+C = 2D+B = 2E+B+C = (n+2) \times B = 360^ \circ$ 
where $n = 1, 2, 3$.

\item[Edge relation:] $a = b = c = d \ne e$.

\item[Heesch number:] $H(T) = 1$.

\item[Corresponding Table and Figure:] Table~\ref{tab05} and Figure~\ref{fig12}.
\end{description}

\begin{table}[!h]
 \begin{center}
{\small
\caption[Table 5]{Value and arrangement of vertices of convex pentagon of Category 5}
\label{tab05}
}
\
{\footnotesize
\begin{tabular}
{c|rrrrr|rrrrr}
\hline
\raisebox{-1.50ex}[0cm][0cm]{$n$}& 
\multicolumn{5}{c|}{\raisebox{-1.75ex}[0.5cm][0.5cm]
{\footnotesize \shortstack{ Value of interior angle \\(degree) }}  } & 
\multicolumn{5}{c}{\raisebox{-1.75ex}[0.5cm][0.5cm]
{\footnotesize \shortstack{ Example of arrangement around each \\vertex (counterclockwise) }} }  \\
 & 

$A$& 
$B$& 
$C$& 
$D$& 
$E$& 
$A$& 
$B$& 
$C$& 
$D$& 
$E$ \\
\hline

1& 
133.47 & 
120 & 
93.07 & 
120 & 
73.47 & 
\textit{ACA}& 
\textit{BBB}& 
\textit{CAA}& 
\textit{DBD}& 
\textit{EECB} \\
\hline

2& 
129.13 & 
90 & 
101.74 & 
135 & 
84.13 & 
\textit{ACA}& 
\textit{BBBB}& 
\textit{CAA}& 
\textit{DBD}& 
\textit{EECB} \\
\hline

3& 
124.74 & 
72 & 
110.51 & 
144 & 
88.74 & 
\textit{ACA}& 
\textit{BBBBB}& 
\textit{CAA}& 
\textit{DBD}& 
\textit{EECB} \\
\hline

\end{tabular}
}
\end{center}
\end{table}

\noindent
\textbf{Remarks.} For $n = 1$, NEEC-spots ``$2B+D = 3D = 2E+C+D = 360^ \circ$'' 
also hold. For $n = 3$, NEEC-spot ``$3B+D = 360^ \circ$'' also holds. The cases of 
$n=0$ and $n \ge 4$ do not exist.

\renewcommand{\figurename}{{\small Figure.}}

\bigskip\bigskip
\bigskip\bigskip
\bigskip\bigskip
\bigskip\bigskip

\bigskip\bigskip
\noindent
\textbf{Category 6}

\begin{figure}[htbp]
 \centering\includegraphics[width=14.5cm,clip]{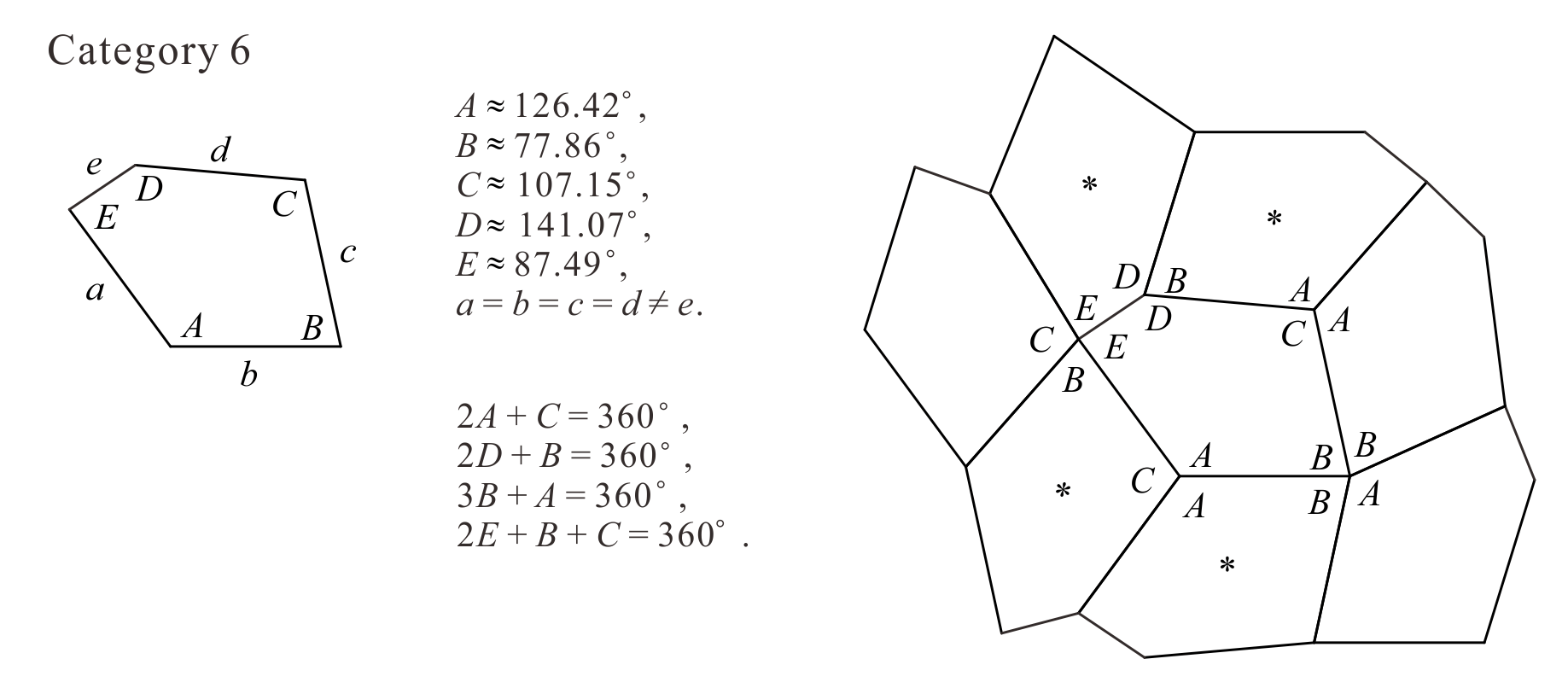} 
  \caption{{\small 
Convex pentagon of Category 6 and example of edge-to-edge corona 
by the pentagons.
} 
\label{fig13}
}
\end{figure}

\begin{description}
 \setlength{\itemindent}{-10pt}
 \setlength{\itemsep}{-3pt} 
\item[Angle relation:] $2A+C = 2D+B = 3B+A = 2E+B+C = 360^ \circ$.

\item[Edge relation:] $a = b = c = d \ne e$.

\item[Heesch number:] $H(T) = 1$.

\item[Corresponding Table and Figure:] Table~\ref{tab06} and Figure~\ref{fig13}.
\end{description}

\begin{table}[!h]
 \begin{center}
{\small
\caption[Table 6]{Value and arrangement of vertices of convex pentagon of Category 6}
\label{tab06}
}
\
{\footnotesize
\begin{tabular}
{rrrrr|rrrrr}
\hline
\multicolumn{5}{c|}{\raisebox{-1.75ex}[0.5cm][0.5cm]
{\footnotesize \shortstack{ Value of interior angle \\(degree) }}  } & 
\multicolumn{5}{c}{\raisebox{-1.75ex}[0.5cm][0.5cm]
{\footnotesize \shortstack{ Example of arrangement around each \\vertex (counterclockwise) }} }  \\

$A$& 
$B$& 
$C$& 
$D$& 
$E$& 
$A$& 
$B$& 
$C$& 
$D$& 
$E$ \\
\hline

126.42 & 
77.86 & 
107.15 & 
141.07 & 
87.49 & 
\textit{ACA}& 
\textit{BBAB}& 
\textit{CAA}& 
\textit{DBD}& 
\textit{EECB} \\
\hline

\end{tabular}
}
\end{center}
\end{table}

\renewcommand{\figurename}{{\small Figure.}}

\noindent
\textbf{Category 7}

\begin{figure}[htbp]
 \centering\includegraphics[width=14.5cm,clip]{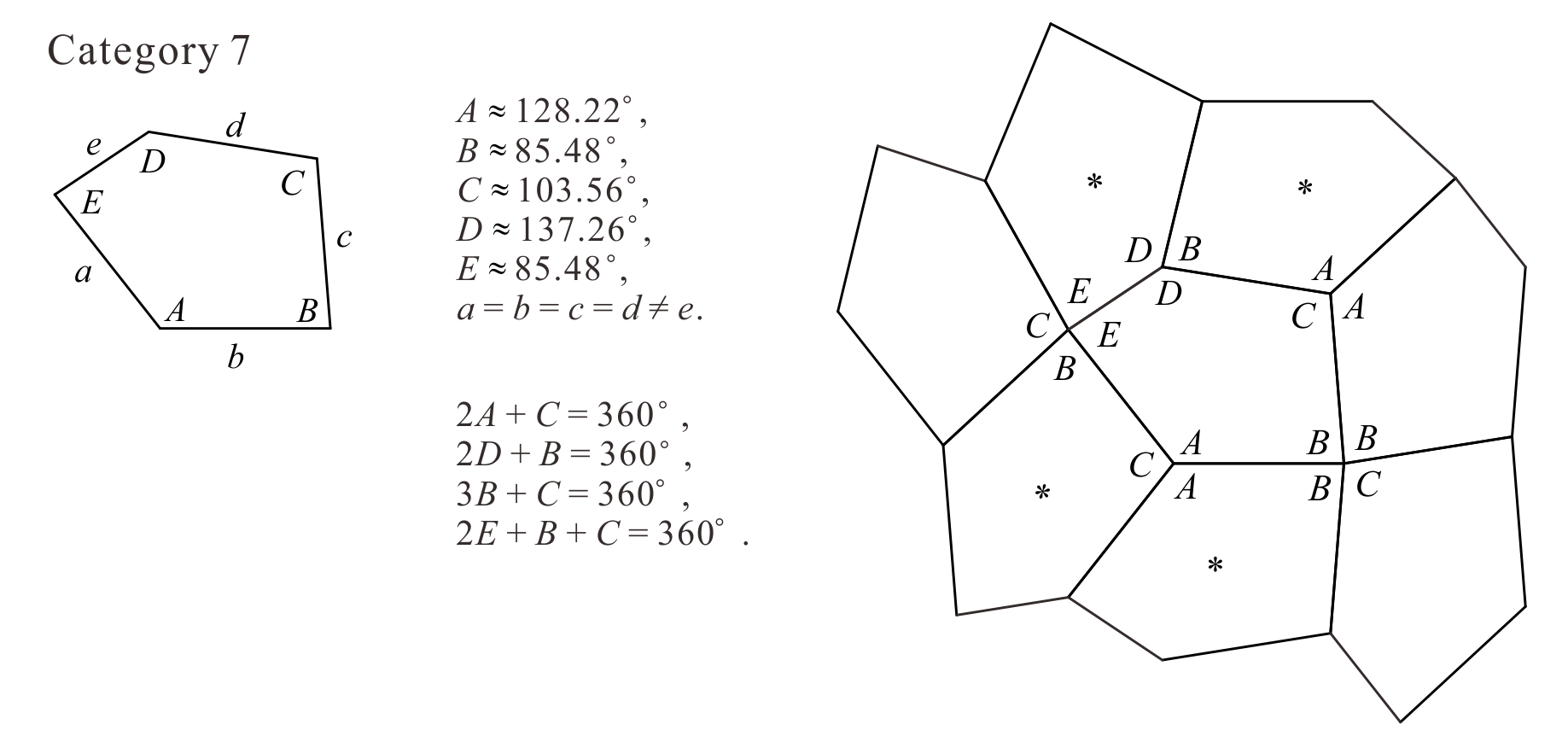} 
  \caption{{\small 
Convex pentagon of Category 7 and example of edge-to-edge corona 
by the pentagons.
} 
\label{fig14}
}
\end{figure}

\begin{description}
 \setlength{\itemindent}{-10pt}
 \setlength{\itemsep}{-3pt} 
\item[Angle relation:] $2A+C = 2D+B = 3B+C = 2E+B+C = 360^ \circ$.

\item[Edge relation:] $a = b = c = d \ne e$.

\item[Heesch number:] $H(T) = 1$.

\item[Corresponding Table and Figure:] Table~\ref{tab07} and Figure~\ref{fig14}.
\end{description}

\begin{table}[!h]
 \begin{center}
{\small
\caption[Table 7]{Value and arrangement of vertices of convex pentagon of Category 7}
\label{tab07}
}
\
{\footnotesize
\begin{tabular}
{rrrrr|rrrrr}
\hline
\multicolumn{5}{c|}{\raisebox{-1.75ex}[0.5cm][0.5cm]
{\footnotesize \shortstack{ Value of interior angle \\(degree) }}  } & 
\multicolumn{5}{c}{\raisebox{-1.75ex}[0.5cm][0.5cm]
{\footnotesize \shortstack{ Example of arrangement around each \\vertex (counterclockwise) }} }  \\

$A$& 
$B$& 
$C$& 
$D$& 
$E$& 
$A$& 
$B$& 
$C$& 
$D$& 
$E$ \\
\hline

128.22 & 
85.48 & 
103.56 & 
137.26 & 
85.48 & 
\textit{ACA}& 
\textit{BBCB}& 
\textit{CAA}& 
\textit{DBD}& 
\textit{EECB} \\
\hline

\end{tabular}
}
\end{center}
\end{table}

\noindent
\textbf{Remarks.} NEEC-spots ``$2D+E=2B+C+E=3E+C=360^ \circ$'' also hold.

\renewcommand{\figurename}{{\small Figure.}}

\bigskip\bigskip
\noindent
\textbf{Category 8}

\begin{figure}[htbp]
 \centering\includegraphics[width=14.5cm,clip]{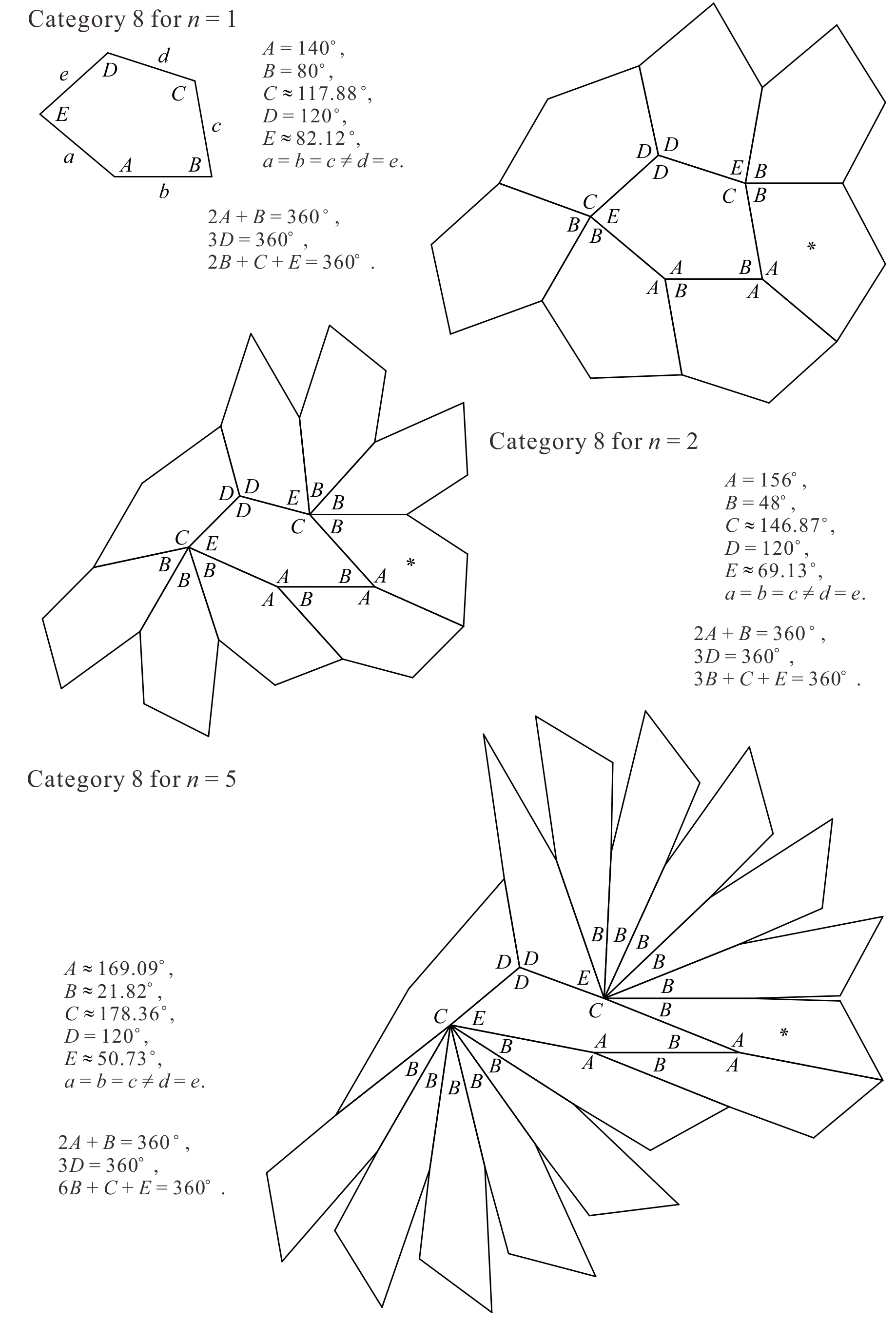} 
  \caption{{\small 
Convex pentagons of Category 8 and examples of edge-to-edge coronas 
by the pentagons.
} 
\label{fig15}
}
\end{figure}

\begin{description}
 \setlength{\itemindent}{-10pt}
 \setlength{\itemsep}{-3pt} 
\item[Angle relation:] $2A+B = 3D = (n+1) \times B+C+E = 360^ \circ$ 
where $n = 1, 2, 3, 4, 5$.

\item[Edge relation:] $a = b = c \ne d = e$.

\item[Heesch number:] $H(T) = 1$.

\item[Corresponding Table and Figure:] Table~\ref{tab08} and Figure~\ref{fig15}.
\end{description}

\begin{table}[!h]
 \begin{center}
{\small
\caption[Table 8]{Value and arrangement of vertices of convex pentagon of Category 8}
\label{tab08}
}
\
{\footnotesize
\begin{tabular}
{c|rrrrr|rrrrr}
\hline
\raisebox{-1.50ex}[0cm][0cm]{$n$}& 
\multicolumn{5}{c|}{\raisebox{-1.75ex}[0.5cm][0.5cm]
{\footnotesize \shortstack{ Value of interior angle \\(degree) }}  } & 
\multicolumn{5}{c}{\raisebox{-1.75ex}[0.5cm][0.5cm]
{\footnotesize \shortstack{ Example of arrangement around each \\vertex (counterclockwise) }} }  \\
 & 

$A$& 
$B$& 
$C$& 
$D$& 
$E$& 
$A$& 
$B$& 
$C$& 
$D$& 
$E$ \\
\hline

1& 
140 & 
80 & 
117.88 & 
120 & 
82.12 & 
\textit{AAB}& 
\textit{BAA}& 
\textit{CBBE}& 
\textit{DDD}& 
\textit{ECBB} \\
\hline

2& 
156 & 
48 & 
146.87 & 
120 & 
69.13 & 
\textit{AAB}& 
\textit{BAA}& 
\textit{CBBBE}& 
\textit{DDD}& 
\textit{ECBBB} \\
\hline

3& 
162.86 & 
34.29 & 
162.34 & 
120 & 
60.52 & 
\textit{AAB}& 
\textit{BAA}& 
\textit{CBBBBE}& 
\textit{DDD}& 
\textit{ECBBBB} \\
\hline

4& 
166.67 & 
26.67 & 
171.91 & 
120 & 
54.76 & 
\textit{AAB}& 
\textit{BAA}& 
\textit{CBBBBBE}& 
\textit{DDD}& 
\textit{ECBBBBB} \\
\hline

5& 
169.09 & 
21.82 & 
178.36 & 
120 & 
50.73 & 
\textit{AAB}& 
\textit{BAA}& 
\textit{CBBBBBBE}& 
\textit{DDD}& 
\textit{ECBBBBBB} \\
\hline

\end{tabular}
}
\end{center}
\end{table}

\noindent
\textbf{Remarks. }The cases of $n = 0$ and $n \ge 6$ do not exist. Note 
that, in all cases, NEEC-spot ``$(2n+1) \times B+D=360^ \circ$'' also holds.

\renewcommand{\figurename}{{\small Figure.}}

\bigskip\bigskip
\noindent
\textbf{Category 9}

\begin{figure}[htbp]
 \centering\includegraphics[width=14.5cm,clip]{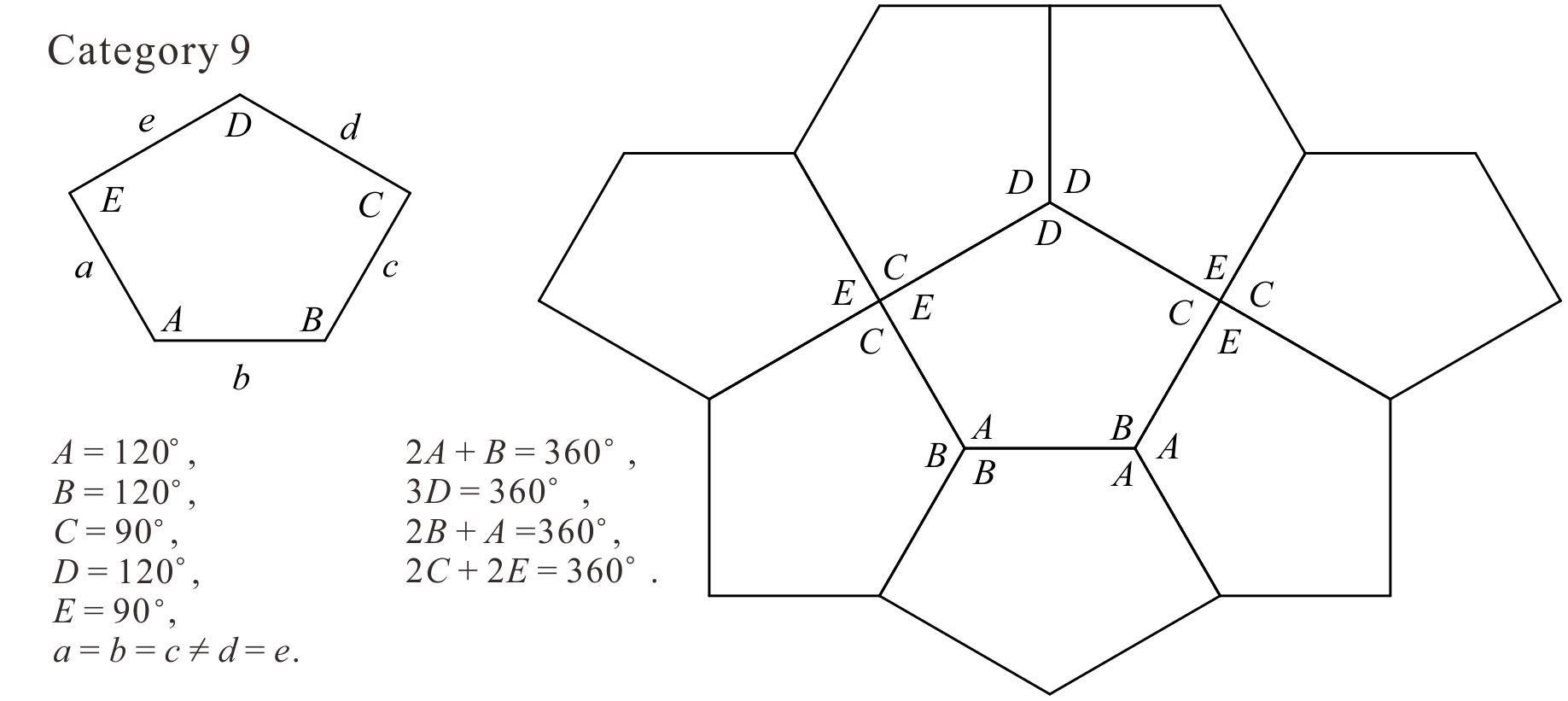} 
  \caption{{\small 
Convex pentagon of Category 9 and example of edge-to-edge corona 
by the pentagons.
} 
\label{fig16}
}
\end{figure}

\begin{description}
 \setlength{\itemindent}{-10pt}
 \setlength{\itemsep}{-3pt} 
\item[Angle relation:] $2A+B = 3D = 2B+A = 2C+2E = 360^ \circ$.

\item[Edge relation:] $a = b = c \ne d = e$.

\item[Heesch number:] $H(T) = 1$.

\item[Corresponding Table and Figure:] Table~\ref{tab09} and Figure~\ref{fig16}.
\end{description}

\begin{table}[!h]
 \begin{center}
{\small
\caption[Table 9]{Value and arrangement of vertices of convex pentagon of Category 9}
\label{tab09}
}
\
{\footnotesize
\begin{tabular}
{rrrrr|rrrrr}
\hline
\multicolumn{5}{c|}{\raisebox{-1.75ex}[0.5cm][0.5cm]
{\small \shortstack{ Value of interior angle \\(degree) }}  } & 
\multicolumn{5}{c}{\raisebox{-1.75ex}[0.5cm][0.5cm]
{\small \shortstack{ Example of arrangement around each \\vertex (counterclockwise) }} }  \\

$A$& 
$B$& 
$C$& 
$D$& 
$E$& 
$A$& 
$B$& 
$C$& 
$D$& 
$E$ \\
\hline

120 & 
120 & 
90 & 
120 & 
90 & 
\textit{ABB}& 
\textit{BAA}& 
\textit{CECE}& 
\textit{DDD}& 
\textit{ECEC} \\
\hline

\end{tabular}
}
\end{center}
\end{table}

\noindent
\textbf{Remarks.} The convex pentagon is symmetric to itself with respect to 
the line through the vertex $D$ and the midpoint of the edge $b$ (\textit{AB}). When the 
vertices $D$ of three convex pentagons are concentrated at a point, it is the 
convex hexagon that cannot generate a tiling~\cite{Sugimoto_2012}. Then, EEC-spots 
``$3A = 3B = 4C = 3C+E = 3E+C = 4E = 360^ \circ$'' and NEEC-spots 
``$A+B+D = 2A+D = 2D+A = 2B+D = 2D+B = 360^ \circ$'' 
also hold.

\renewcommand{\figurename}{{\small Figure.}}

\bigskip\bigskip
\noindent
\textbf{Category 10}

\begin{figure}[htbp]
 \centering\includegraphics[width=14.5cm,clip]{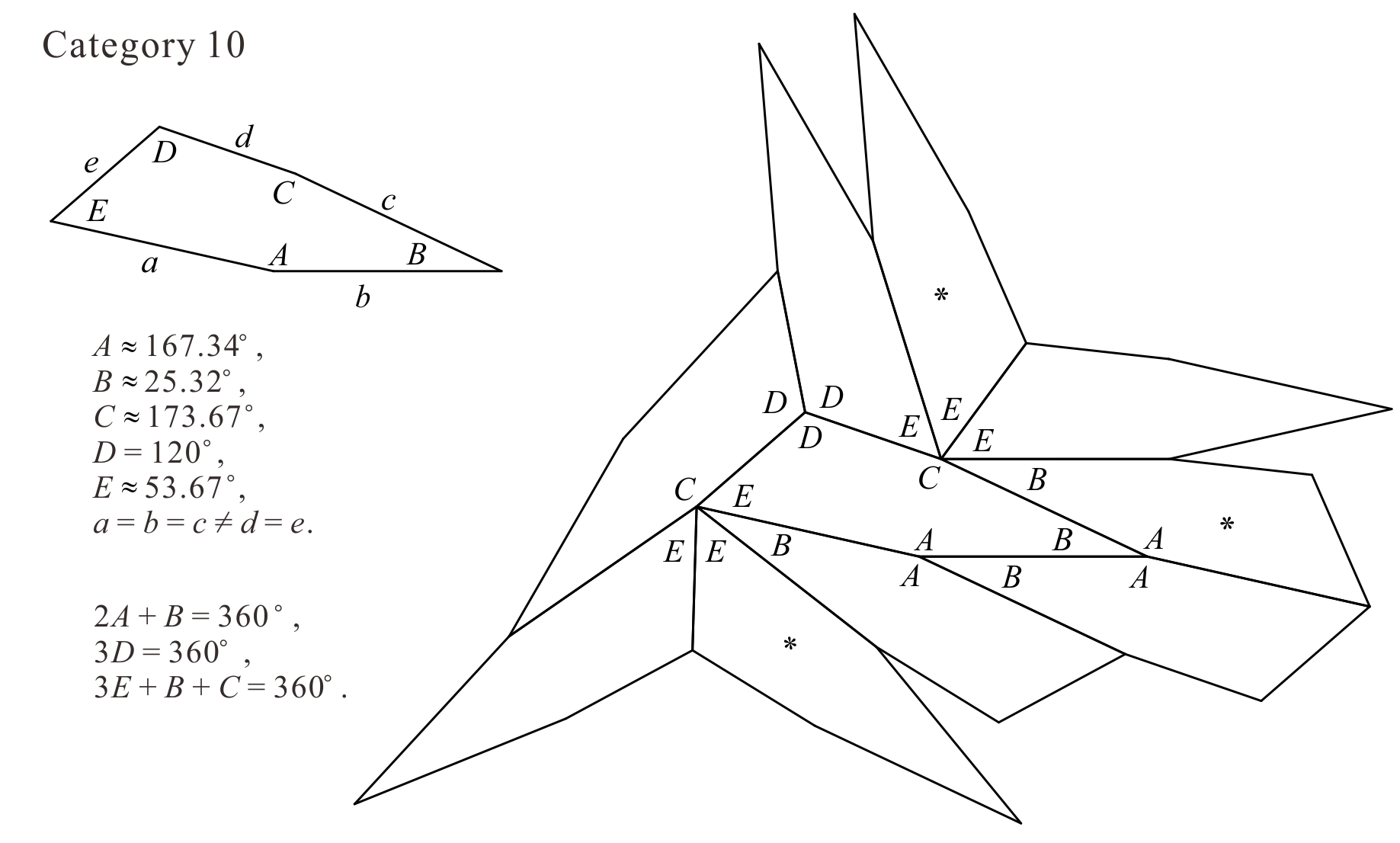} 
  \caption{{\small 
Convex pentagon of Category 10 and example of edge-to-edge corona 
by the pentagons.
} 
\label{fig17}
}
\end{figure}

\begin{description}
 \setlength{\itemindent}{-10pt}
 \setlength{\itemsep}{-3pt} 
\item[Angle relation:] $2A+B = 3D = 3E+B+C = 360^ \circ$.

\item[Edge relation:] $a = b = c \ne d = e$.

\item[Heesch number:] $H(T) = 1$.

\item[Corresponding Table and Figure:] Table~\ref{tab10} and Figure~\ref{fig17}.
\end{description}

\begin{table}[!h]
 \begin{center}
{\small
\caption[Table 10]{Value and arrangement of vertices of convex pentagon of Category 10}
\label{tab10}
}
\
{\footnotesize
\begin{tabular}
{rrrrr|rrrrr}
\hline
\multicolumn{5}{c|}{\raisebox{-1.75ex}[0.5cm][0.5cm]
{\small \shortstack{ Value of interior angle \\(degree) }}  } & 
\multicolumn{5}{c}{\raisebox{-1.75ex}[0.5cm][0.5cm]
{\small \shortstack{ Example of arrangement around each \\vertex (counterclockwise) }} }  \\

$A$& 
$B$& 
$C$& 
$D$& 
$E$& 
$A$& 
$B$& 
$C$& 
$D$& 
$E$ \\
\hline

167.34 & 
25.32 & 
173.67 & 
120 & 
53.67 & 
\textit{AAB}& 
\textit{BAA}& 
\textit{CBEEE}& 
\textit{DDD}& 
\textit{ECEEB} \\
\hline

\end{tabular}
}
\end{center}
\end{table}

\noindent
\textbf{Remarks.} EEC-spot ``$4E+B+D=360^ \circ$'' also holds.

\renewcommand{\figurename}{{\small Figure.}}

\bigskip\bigskip
\noindent
\textbf{Category 11}

\begin{figure}[htbp]
 \centering\includegraphics[width=14.5cm,clip]{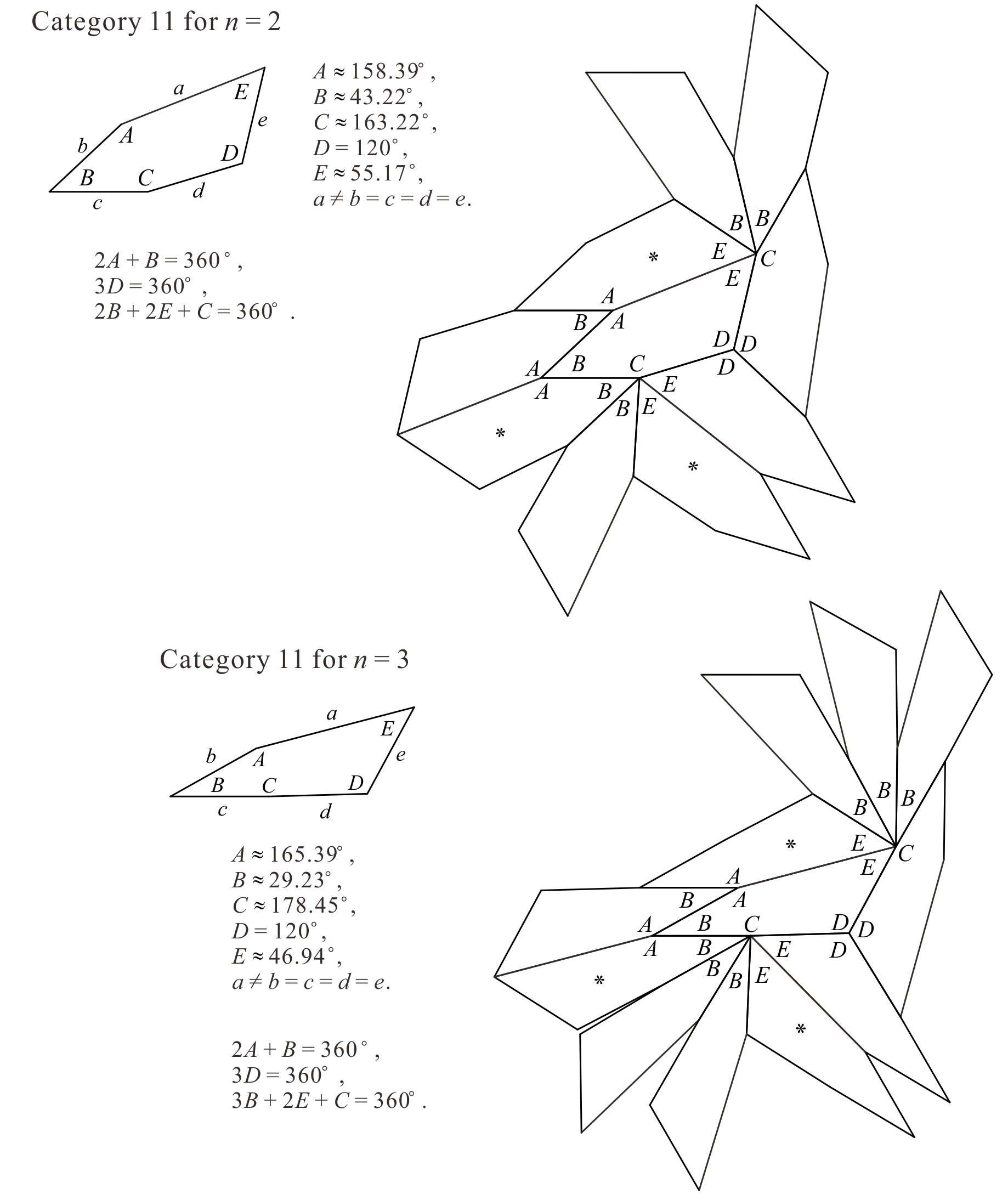} 
  \caption{{\small 
Convex pentagons of Category 11 and examples of edge-to-edge coronas 
by the pentagons.
} 
\label{fig18}
}
\end{figure}

\begin{description}
 \setlength{\itemindent}{-10pt}
 \setlength{\itemsep}{-3pt} 
\item[Angle relation:] $2A+B = 3D= n \times B+2E+C = 360^ \circ$ where $n = 1, 2, 3$.

\item[Edge relation:] $a \ne b = c = d = e$.

\item[Heesch number:]  $H(T) = \infty$ for $n=1$ or $H(T) = 1$ for $n = 2, 3$.

\item[Corresponding Table and Figure:] Table~\ref{tab11} and Figure~\ref{fig18}.
\end{description}

\begin{table}[!h]
 \begin{center}
{\small
\caption[Table 11]{Value and arrangement of vertices of convex pentagon of Category 11}
\label{tab11}
}
\
{\footnotesize
\begin{tabular}
{c|c|rrrrr|rrrrr}
\hline
\raisebox{-1.50ex}[0cm][0cm]{$n$}& 
\raisebox{-1.50ex}[0cm][0cm]{$H(T)$}& 
\multicolumn{5}{c|}{\raisebox{-1.75ex}[0.5cm][0.5cm]
{\footnotesize \shortstack{ Value of interior angle \\(degree) }}  } & 
\multicolumn{5}{c}{\raisebox{-1.75ex}[0.5cm][0.5cm]
{\footnotesize \shortstack{ Example of arrangement around each \\vertex (counterclockwise) }} }  \\
 & 
 & 

$A$& 
$B$& 
$C$& 
$D$& 
$E$& 
$A$& 
$B$& 
$C$& 
$D$& 
$E$ \\
\hline
 
1& 
$\infty $& 
139.11 & 
81.78 & 
120 & 
120 & 
79.11 & 
\multicolumn{5}{c}{Convex Pentagonal tile belonging to Type 8}  \\

\hline

2& 
1& 
158.39 & 
43.22 & 
163.22 & 
120 & 
55.17 & 
\textit{AAB}& 
\textit{BAA}& 
\textit{CBBEE}& 
\textit{DDD}& 
\textit{ECBBE} \\
\hline

3& 
1& 
165.39 & 
29.23 & 
178.45 & 
120 & 
46.94 & 
\textit{AAB}& 
\textit{BAA}& 
\textit{CBBBEE}& 
\textit{DDD}& 
\textit{ECBBBE} \\
\hline

\end{tabular}
}
\end{center}
\end{table}

\noindent
\textbf{Remarks.} As for the convex pentagon satisfying 
``$2A + B = 3D = 360^ \circ , a \ne b = c = d = e$,'' if the convex pentagon 
also has the relation ``$2E+B+C = 360^ \circ$,'' then it is the convex pentagonal 
tile belonging to Type 8 (in this case, EEC-spots 
``$3C = 2C+D = 2D+C = 2E+B+D = 360^ \circ$'' also hold). 
Therefore, the convex pentagon of Category 11 is the convex pentagon with 
$H(T) = \infty $ if and only if $n = 1$. Note that, in all cases, EEC-spot 
``$(2n - 1)\times B + 2E + D = 360^ \circ $'' should also holds. The cases 
of $n=0$ and $n \ge 4$ do not exist.

\renewcommand{\figurename}{{\small Figure.}}

\bigskip\bigskip
\bigskip\bigskip

\bigskip\bigskip
\noindent
\textbf{Category 12}

\begin{figure}[htbp]
 \centering\includegraphics[width=13.2cm,clip]{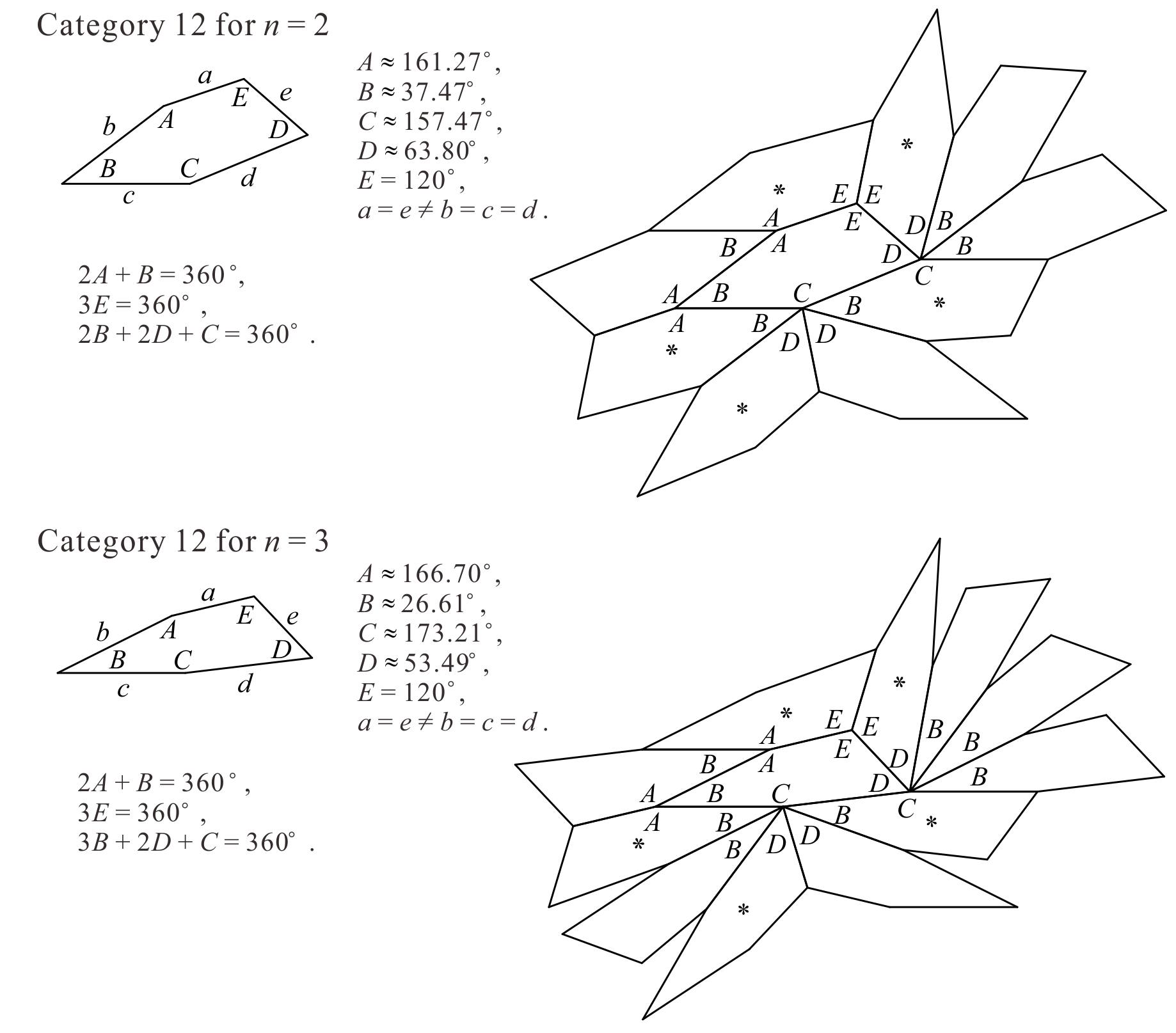} 
  \caption{{\small 
Convex pentagons of Category 12 and examples of edge-to-edge coronas 
by the pentagons.
} 
\label{fig19}
}
\end{figure}

\begin{description}
 \setlength{\itemindent}{-10pt}
 \setlength{\itemsep}{-3pt} 
\item[Angle relation:] $2A+B = 3E= n \times B+2D+C = 360$ where $n = 1, 2, 3$.

\item[Edge relation:] $a = e \ne b = c = d$.

\item[Heesch number:] $H(T) = \infty$ for $n=1$ or $H(T) = 1$ for $n = 2, 3$.

\item[Corresponding Table and Figure:] Table~\ref{tab12} and Figure~\ref{fig19}.
\end{description}

\begin{table}[!h]
 \begin{center}
{\small
\caption[Table 12]{Value and arrangement of vertices of convex pentagon of Category 12}
\label{tab12}
}
\
{\footnotesize
\begin{tabular}
{c|c|rrrrr|rrrrr}
\hline
\raisebox{-1.50ex}[0cm][0cm]{$n$}& 
\raisebox{-1.50ex}[0cm][0cm]{$H(T)$}& 
\multicolumn{5}{c|}{\raisebox{-1.75ex}[0.5cm][0.5cm]
{\footnotesize \shortstack{ Value of interior angle \\(degree) }}  } & 
\multicolumn{5}{c}{\raisebox{-1.75ex}[0.5cm][0.5cm]
{\footnotesize \shortstack{ Example of arrangement around each \\vertex (counterclockwise) }} }  \\
 & 
 & 

$A$& 
$B$& 
$C$& 
$D$& 
$E$& 
$A$& 
$B$& 
$C$& 
$D$& 
$E$ \\
\hline
 
1& 
$\infty $& 
150 & 
60 & 
120 & 
90 & 
120 &  
\multicolumn{5}{c}{Convex Pentagonal tile belonging to Types 1 and 5}  \\

\hline

2& 
1& 
161.27 & 
37.47 & 
157.47 & 
63.80 & 
120 & 
\textit{AAB}& 
\textit{BAA}& 
\textit{CBDDB}& 
\textit{DCBBD}& 
\textit{EEE} \\
\hline

3& 
1& 
166.70 & 
26.61 & 
173.21 & 
53.49 & 
120 & 
\textit{AAB}& 
\textit{BAA}& 
\textit{CBBDDB}& 
\textit{DCBBBD}& 
\textit{EEE} \\
\hline

\end{tabular}
}
\end{center}
\end{table}

\noindent
\textbf{Remarks.} As for the convex pentagon satisfying 
``$2A + B = 3E = 360^ \circ ,a = e \ne b = c = d$,'' if the convex pentagon 
also has the relation ``$ 2D+B+C = 360^ \circ$,'' then it is the convex 
pentagonal tile belonging to Types 1 and 5 (in this case, EEC-spots 
``$A+C+D =A+D+E = 3C = 2B+A+D = 2B+2C = 2D+B+E = 4D = 
4B+C = 3B+2D = 6B = 360^ \circ$'' and NEE-spots 
``$2C+E = 2B+C+E = 2B+2E = 4B+E = 360^ \circ$'' also hold). That is, 
the convex pentagon of Category 12 is the convex pentagon 
with $H(T) = \infty $ if and only if $n = 1$. We called this convex pentagon 
of $n = 1$ the ``TH-pentagon,'' and it has the relation of 
heptiamonds~\cite{Su_and_Ar_2017}. 
Note that, in all cases, EEC-spot ``$2(n - 1)\times B + 2D + E = 360^ \circ $'' 
should also holds. The cases of $n=0$ and $n \ge 4$ do not exist.

\renewcommand{\figurename}{{\small Figure.}}

\renewcommand{\figurename}{{\small Figure.}}
\begin{figure}[htbp]
 \centering\includegraphics[width=13.2cm,clip]{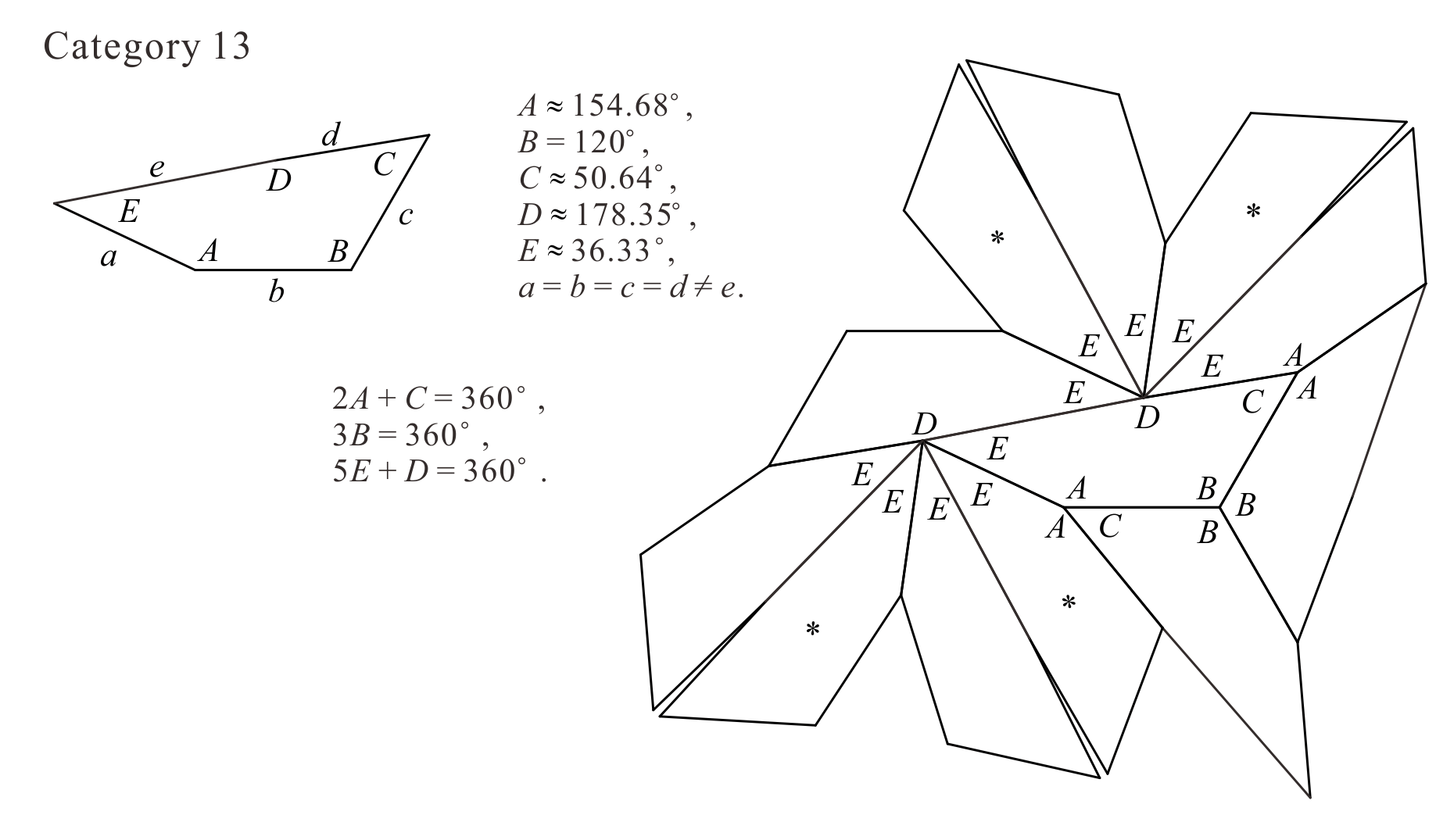} 
  \caption{{\small 
Convex pentagon of Category 13 and example of edge-to-edge corona 
by the pentagons.
} 
\label{fig20}
}
\end{figure}

\bigskip\bigskip
\noindent
\textbf{Category 13}

\begin{description}
 \setlength{\itemindent}{-10pt}
 \setlength{\itemsep}{-3pt} 
\item[Angle relation:] $2A+C = 3B = 5E+D = 360^ \circ$.

\item[Edge relation:] $a = b = c = d \ne e$.

\item[Heesch number:] $H(T) = 1$.

\item[Corresponding Table and Figure:] Table~\ref{tab13} and Figure~\ref{fig20}.
\end{description}

\begin{table}[!h]
 \begin{center}
{\small
\caption[Table 13]{Value and arrangement of vertices of convex pentagon of Category 13}
\label{tab13}
}
\
{\footnotesize
\begin{tabular}
{rrrrr|rrrrr}
\hline
\multicolumn{5}{c|}{\raisebox{-1.75ex}[0.5cm][0.5cm]
{\small \shortstack{ Value of interior angle \\(degree) }}  } & 
\multicolumn{5}{c}{\raisebox{-1.75ex}[0.5cm][0.5cm]
{\small \shortstack{ Example of arrangement around each \\vertex (counterclockwise) }} }  \\

$A$& 
$B$& 
$C$& 
$D$& 
$E$& 
$A$& 
$B$& 
$C$& 
$D$& 
$E$ \\
\hline

154.68 & 
120 & 
50.64 & 
178.35 & 
36.33 & 
\textit{AAC}& 
\textit{BBB}& 
\textit{CAA}& 
\textit{DEEEEE}& 
\textit{EDEEEE} \\
\hline

\end{tabular}
}
\end{center}
\end{table}

\bigskip\bigskip
\bigskip\bigskip
\bigskip\bigskip
\bigskip\bigskip
\bigskip\bigskip

\bigskip\bigskip
\noindent
\textbf{Category 14}

\begin{figure}[htbp]
 \centering\includegraphics[width=14.5cm,clip]{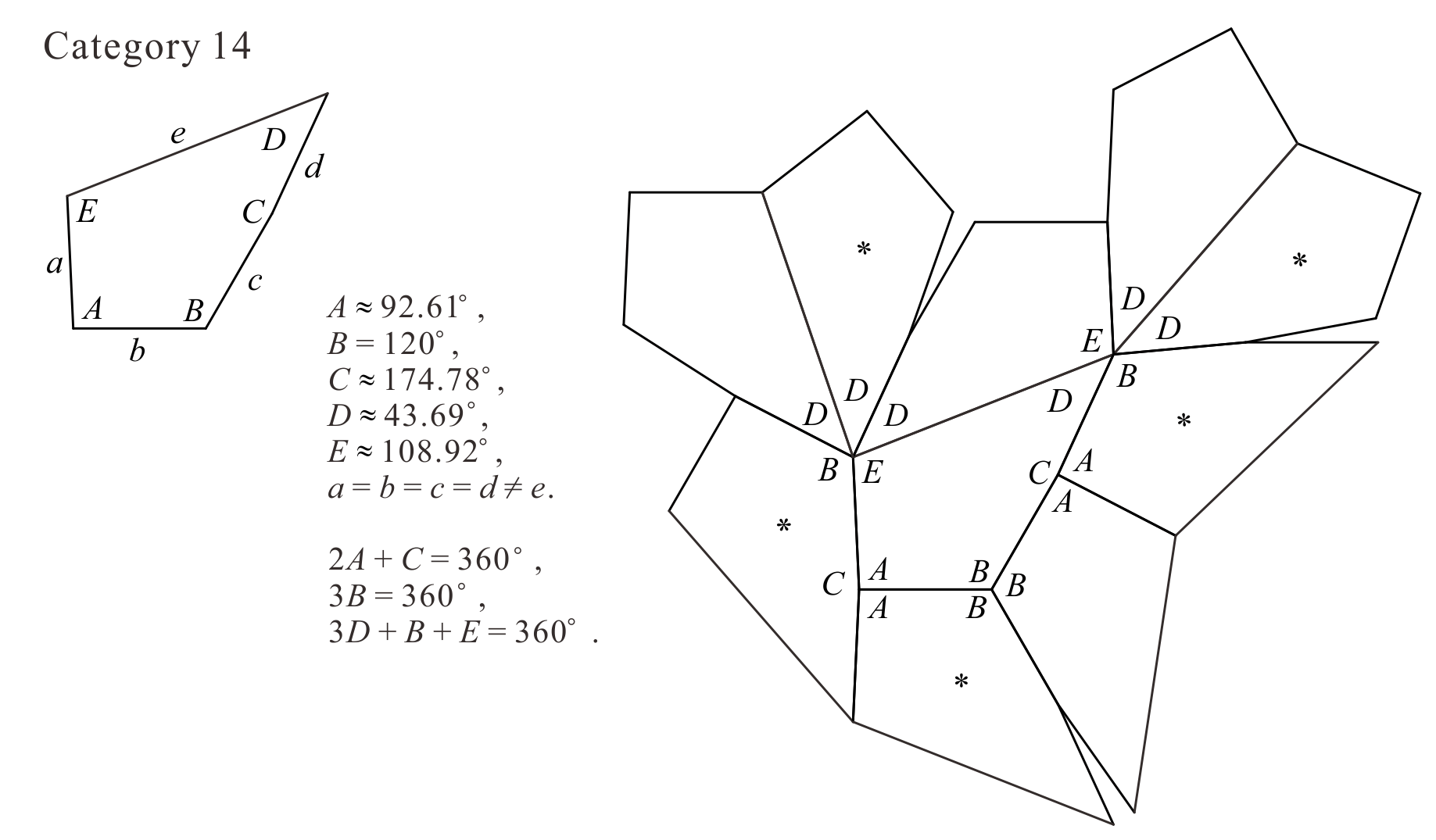} 
  \caption{{\small 
Convex pentagon of Category 14 and example of edge-to-edge corona 
by the pentagons.
} 
\label{fig21}
}
\end{figure}

\begin{description}
 \setlength{\itemindent}{-10pt}
 \setlength{\itemsep}{-3pt} 
\item[Angle relation:] $2A+C = 3B = 3D+B+E = 360^ \circ$.

\item[Edge relation:] $a = b = c = d \ne e$.

\item[Heesch number:] $H(T) = 1$.

\item[Corresponding Table and Figure:] Table~\ref{tab14} and Figure~\ref{fig21}.
\end{description}

\begin{table}[!h]
 \begin{center}
{\small
\caption[Table 14]{Value and arrangement of vertices of convex pentagon of Category 14}
\label{tab14}
}
\
{\footnotesize
\begin{tabular}
{rrrrr|rrrrr}
\hline
\multicolumn{5}{c|}{\raisebox{-1.75ex}[0.5cm][0.5cm]
{\small \shortstack{ Value of interior angle \\(degree) }}  } & 
\multicolumn{5}{c}{\raisebox{-1.75ex}[0.5cm][0.5cm]
{\small \shortstack{ Example of arrangement around each \\vertex (counterclockwise) }} }  \\

$A$& 
$B$& 
$C$& 
$D$& 
$E$& 
$A$& 
$B$& 
$C$& 
$D$& 
$E$ \\
\hline

92.61 & 
120 & 
174.78 & 
43.69 & 
108.92 & 
\textit{ACA}& 
\textit{BBB}& 
\textit{CAA}& 
\textit{DBDDE}& 
\textit{EDDDB} \\
\hline

\end{tabular}
}
\end{center}
\end{table}

\noindent
\textbf{Remarks}. NEEC-spot ``$4D+2A=360^ \circ$'' also holds.

\renewcommand{\figurename}{{\small Figure.}}

\bigskip\bigskip
\noindent
\textbf{Category 15}

\begin{figure}[htb]
 \centering\includegraphics[width=14.5cm,clip]{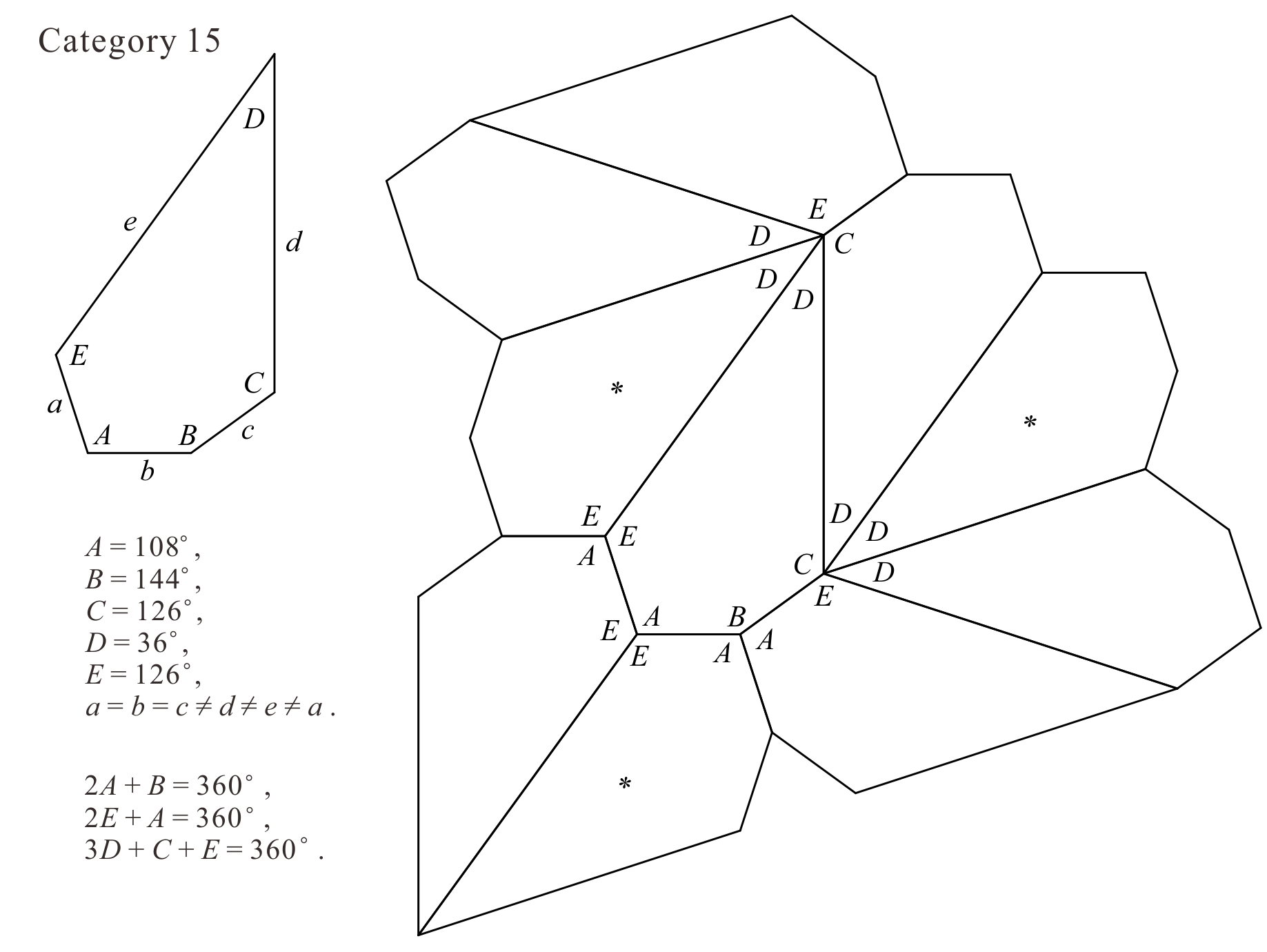} 
  \caption{{\small 
Convex pentagon of Category 15 and example of edge-to-edge corona 
by the pentagons.
} 
\label{fig22}
}
\end{figure}

\begin{description}
 \setlength{\itemindent}{-10pt}
 \setlength{\itemsep}{-3pt} 
\item[Angle relation:] $2A+B = 2E+A = 3D+C+E = 360^ \circ$.

\item[Edge relation:] $a = b = c \ne d \ne e \ne a$.

\item[Heesch number:] $H(T) = 1$.

\item[Corresponding Table and Figure:] Table~\ref{tab15} and Figure~\ref{fig22}.
\end{description}

\begin{table}[!h]
 \begin{center}
{\small
\caption[Table 15]{Value and arrangement of vertices of convex pentagon of Category 15}
\label{tab15}
}
\
{\footnotesize
\begin{tabular}
{rrrrr|rrrrr}
\hline
\multicolumn{5}{c|}{\raisebox{-1.75ex}[0.5cm][0.5cm]
{\small \shortstack{ Value of interior angle \\(degree) }}  } & 
\multicolumn{5}{c}{\raisebox{-1.75ex}[0.5cm][0.5cm]
{\small \shortstack{ Example of arrangement around each \\vertex (counterclockwise) }} }  \\

$A$& 
$B$& 
$C$& 
$D$& 
$E$& 
$A$& 
$B$& 
$C$& 
$D$& 
$E$ \\
\hline

108 & 
144 & 
126 & 
36 & 
126 & 
\textit{AEE}& 
\textit{BAA}& 
\textit{CEDDD}& 
\textit{DCEDD}& 
\textit{EEA} \\
\hline

\end{tabular}
}
\end{center}
\end{table}

\bigskip\bigskip
\bigskip\bigskip
\bigskip\bigskip
\bigskip\bigskip
\bigskip\bigskip
\bigskip\bigskip
\bigskip\bigskip
\bigskip\bigskip

\noindent
\textbf{Remarks.} EEC-spot ``$2C+A = 10D = 360^ \circ$'' and NEEC-spots 
``$A+C+E = 3A+D = 2B+2D = 3D+A+B = 3D+2C = 3D+2E = 4D+2A = 
6D+B = 7D+A=360^ \circ$'' also hold.

\renewcommand{\figurename}{{\small Figure.}}

\bigskip\bigskip
\noindent
\textbf{Category 16}

\begin{figure}[htbp]
 \centering\includegraphics[width=14.5cm,clip]{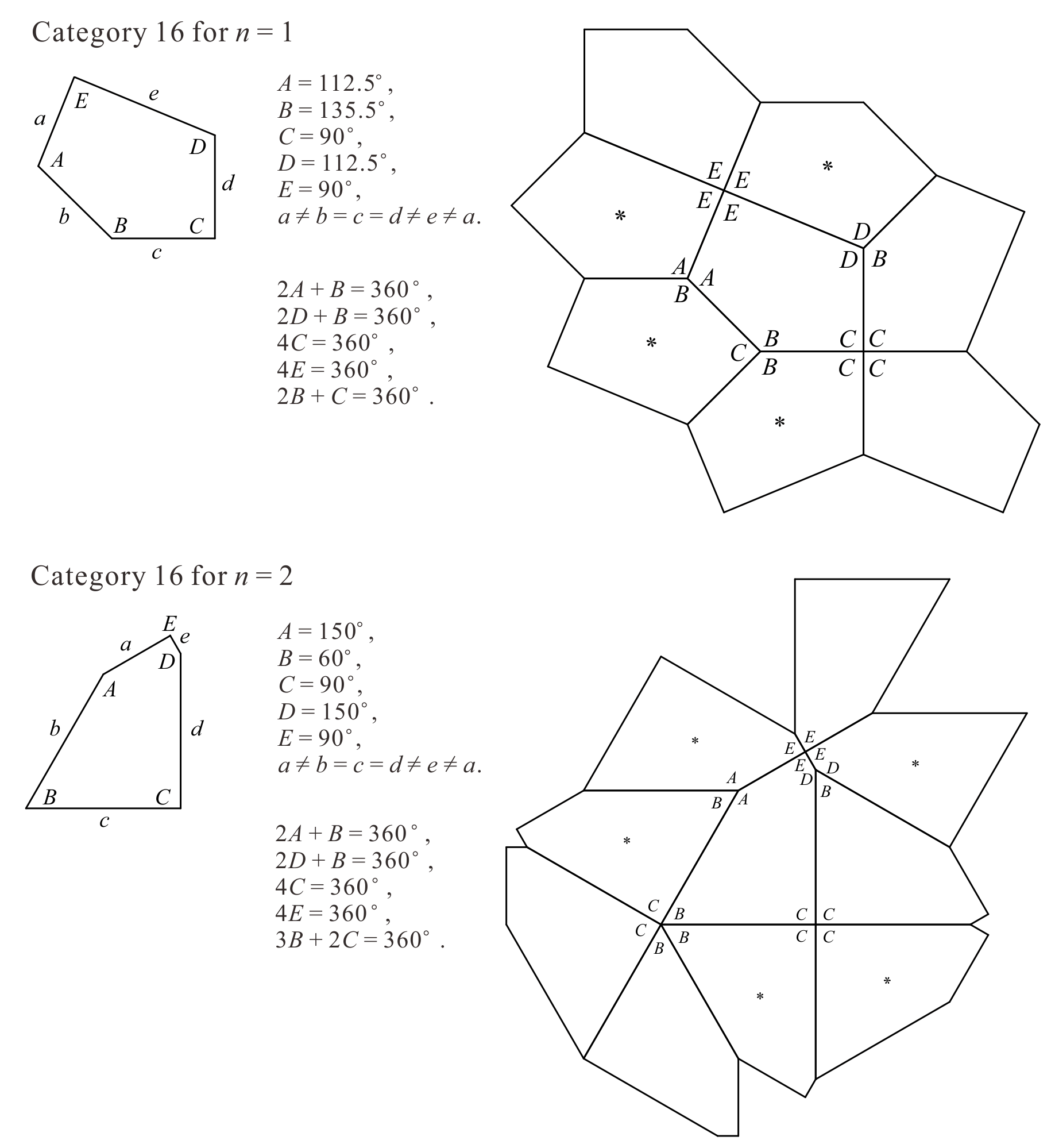} 
  \caption{{\small 
Convex pentagons of Category 16 and examples of edge-to-edge coronas 
by the pentagons.
} 
\label{fig23}
}
\end{figure}

\begin{description}
 \setlength{\itemindent}{-10pt}
 \setlength{\itemsep}{-3pt} 
\item[Angle relation:] $2A+B = 2D+B = 4C = 4E = (n+1) \times B+n \times C = 360^ \circ$ 
where $n = 1, 2$.

\item[Edge relation:] $a \ne b = c = d \ne e \ne a$.

\item[Heesch number:] $H(T) = 1$.

\item[Corresponding Table and Figure:] Table~\ref{tab16} and Figure~\ref{fig23}.
\end{description}

\begin{table}[!h]
 \begin{center}
{\small
\caption[Table 16]{Value and arrangement of vertices of convex pentagon of Category 16}
\label{tab16}
}
\
{\footnotesize
\begin{tabular}
{c|rrrrr|rrrrr}
\hline
\raisebox{-1.50ex}[0cm][0cm]{$n$}& 
\multicolumn{5}{c|}{\raisebox{-1.75ex}[0.5cm][0.5cm]
{\footnotesize \shortstack{ Value of interior angle \\(degree) }}  } & 
\multicolumn{5}{c}{\raisebox{-1.75ex}[0.5cm][0.5cm]
{\footnotesize \shortstack{ Example of arrangement around each \\vertex (counterclockwise) }} }  \\
 & 

$A$& 
$B$& 
$C$& 
$D$& 
$E$& 
$A$& 
$B$& 
$C$& 
$D$& 
$E$ \\
\hline

1& 
112.5 & 
135 & 
90 & 
112.5 & 
90 & 
\textit{AAB}& 
\textit{BCB}& 
\textit{CCCC}& 
\textit{DBD}& 
\textit{EEEE} \\
\hline

2& 
150 & 
60 & 
90 & 
150 & 
90 & 
\textit{AAB}& 
\textit{BCCBB}& 
\textit{CCCC}& 
\textit{DBD}& 
\textit{EEEE} \\
\hline

\end{tabular}
}
\end{center}
\end{table}

\bigskip\bigskip
\bigskip\bigskip
\bigskip\bigskip
\bigskip\bigskip
\bigskip\bigskip
\bigskip\bigskip
\bigskip\bigskip
\bigskip\bigskip

\noindent
\textbf{Remarks.} For $n = 1$, NEEC-spots 
``$A+B+D = 2B+E = 3C+E = 2C+2E = 3E+C = 360^ \circ$'' also hold. 
For $n = 2$, EEC-spot ``$6B = 360^ \circ$'' and NEEC-spots 
``$A+B+D = 2B+A+C = 2B+A+E = 2B+C+D = 2B+D+E = 3C+E = 
2C+2E = 3E+C = 3B+C+E = 3B+2E = 360^ \circ$'' also hold. The cases 
of $n=0$ and o$n \ge 3$ do not exist.

\renewcommand{\figurename}{{\small Figure.}}

\bigskip\bigskip
\noindent
\textbf{Category 17}

\begin{figure}[htbp]
 \centering\includegraphics[width=14.5cm,clip]{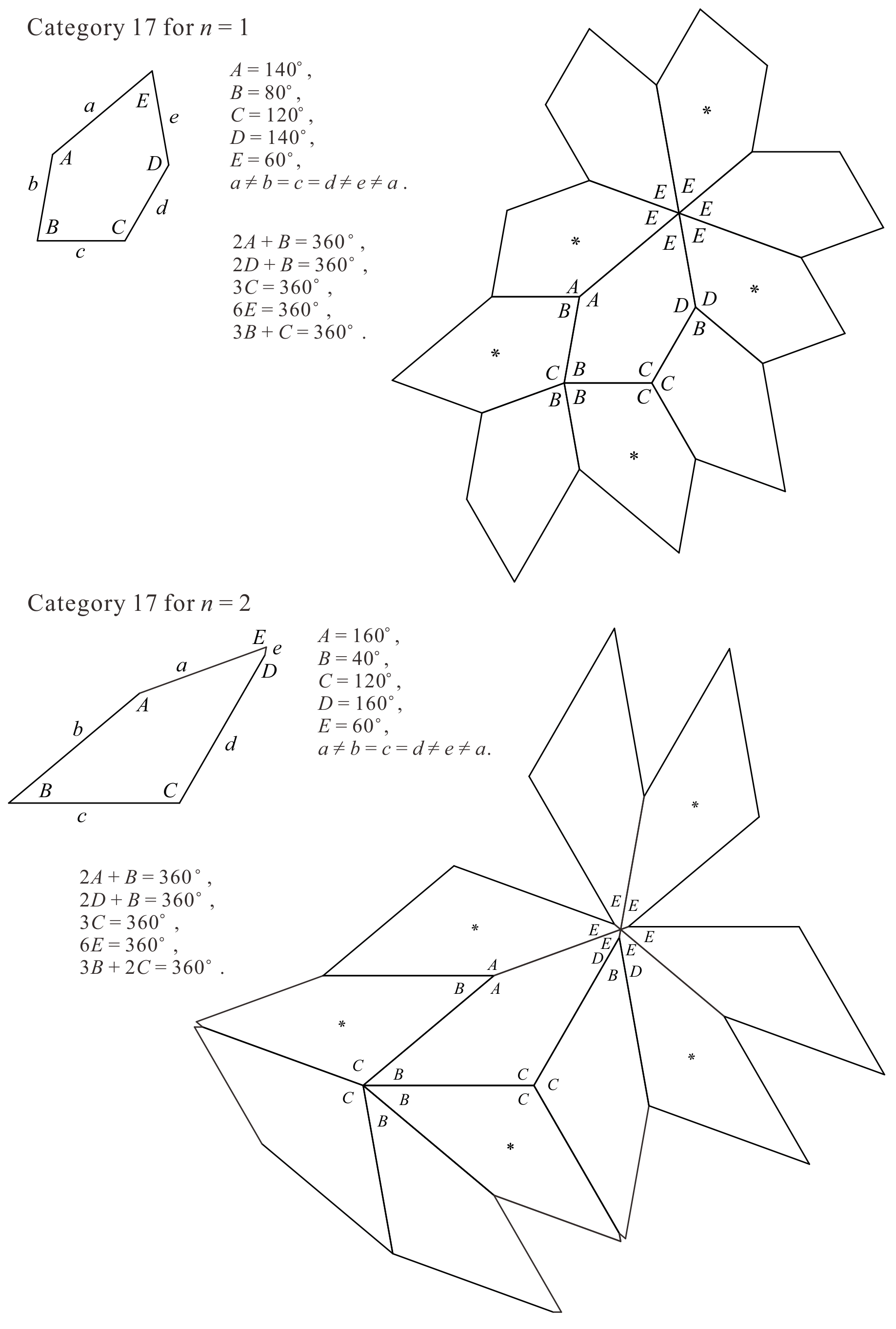} 
  \caption{{\small 
Convex pentagons of Category 17 and examples of edge-to-edge coronas 
by the pentagons.
} 
\label{fig24}
}
\end{figure}

\begin{description}
 \setlength{\itemindent}{-10pt}
 \setlength{\itemsep}{-3pt} 
\item[Angle relation:] $2A+B = 2D+B = 3C = 6E = n \times C + 3B = 360^ \circ$ 
where $n = 1, 2$.

\item[Edge relation:] $a \ne b = c = d \ne e \ne a$.

\item[Heesch number:] $H(T) = 1$.

\item[Corresponding Table and Figure:] Table~\ref{tab17} and Figure~\ref{fig24}.
\end{description}

\begin{table}[!h]
 \begin{center}
{\small
\caption[Table 17]{Value and arrangement of vertices of convex pentagon of Category 17}
\label{tab17}
}
\
{\footnotesize
\begin{tabular}
{c|rrrrr|rrrrr}
\hline
\raisebox{-1.50ex}[0cm][0cm]{$n$}& 
\multicolumn{5}{c|}{\raisebox{-1.75ex}[0.5cm][0.5cm]
{\footnotesize \shortstack{ Value of interior angle \\(degree) }}  } & 
\multicolumn{5}{c}{\raisebox{-1.75ex}[0.5cm][0.5cm]
{\footnotesize \shortstack{ Example of arrangement around each \\vertex (counterclockwise) }} }  \\
 & 

$A$& 
$B$& 
$C$& 
$D$& 
$E$& 
$A$& 
$B$& 
$C$& 
$D$& 
$E$ \\
\hline

1& 
140 & 
80 & 
120 & 
140 & 
60 & 
\textit{AAB}& 
\textit{BCBB}& 
\textit{CCC}& 
\textit{DBD}& 
\textit{EEEEEE} \\
\hline

2& 
160 & 
40 & 
120 & 
160 & 
60 & 
\textit{AAB}& 
\textit{BCCBB}& 
\textit{CCC}& 
\textit{DBD}& 
\textit{EEEEEE} \\
\hline

\end{tabular}
}
\end{center}
\end{table}

\noindent
\textbf{Remarks.} For $n = 1$, NEEC-spots ``$A+B+D = 2B+A+E = 2B+D+E = 
2C+2E = 3B+2E = 4E+C = 360^ \circ$'' also hold. For $n = 2$, EEC-spots 
``$6B+C = 9B = 360^ \circ$'' and NEEC-spots ``$A+B+D = 2B+A+C = 2B+C+D = 
2C+2E = 2B+2E+A = 2B+2E+D = 4E+C = 5B+A = 5B+D = 3B+2E+C = 
6B+2E = 360^ \circ$'' also hold. The cases of $n = 0$ and $n \ge 3$ do not exist.

\renewcommand{\figurename}{{\small Figure.}}

\subsection{Supplement}
\label{subsection2.3}

Category 4 contains the convex pentagonal tile belonging to Type 9, and 
Category 11 contains the convex pentagonal tile belonging to Type 8. It 
might seem that there is no Category containing a convex pentagonal tile 
belonging to Type 7 (see Figure~\ref{fig01}), whose conditions are similar to the 
convex pentagonal tiles of Types 8 and 9. In conclusion, such a Category 
does not exist due to the properties of the convex pentagonal tile of Type 
7. That can be deduced as follows. Let $P_{4}$ be a convex pentagon that 
satisfies the conditions ``$2B+C = 2D+A = 360^ \circ, a = b = c = d$'' of Type 7. 
In $P_{4}$, consider two isosceles triangles, \textit{ABE} with base angles 
$\alpha $ and \textit{BCD} with base angles $\beta $. Then the interior 
angles of $P_{4}$ can be expressed as follows:

\begin{equation}
\label{eq8}
\left\{ {\begin{array}{l}
 A = 180^ \circ - 2\alpha , \\ 
 B = 90^ \circ + \beta , \\ 
 C = 180^ \circ - 2\beta , \\ 
 D = 90^ \circ + \alpha , \\ 
 E = \alpha + \beta . \\ 
 \end{array}} \right.
\end{equation}

\noindent
If $P_{4}$ also has the properties of Categories 4 and 11, it can be 
considered that ``$n\times E + A + C = 360^ \circ $'' holds (corresponding 
to Type 7 when $n = 2)$. However, from (\ref{eq8}), the relation of 
``$n\times E = 360^ \circ - A - C = 2\alpha + 2\beta $'' is admitted 
geometrically only in the case of $n = 2$. Therefore, no Category 
contains a convex pentagon tile belonging to Type 7.

As of date, it is unknown whether there is a Category that includes convex 
pentagonal tiles of Types 2 or 4 that can generate edge-to-edge tilings.

\section{Conclusions}
\label{section3}

In this manuscript, we found a convex pentagon with $H(T) = 1$ admitting 
edge-to-edge corona that was classified based on Category. In addition, 
among the convex pentagons with $H(T) = 1$ found in this survey, it is 
possible to form a first corona even without using a convex pentagon, in 
which only the convex pentagon of Category 9 is reversed.

The results of this manuscript were obtained by a survey targeting the 
convex pentagons shown in Table 1 of \cite{Sugimoto_2015}\footnote{ We got 
most of the results of this manuscript in 2015. On the other hand, Craig S. Kaplan 
independently discovered infinite families of convex pentagons corresponding to 
Category 3 in 2017~\cite{Kaplan_2017_p4}.}. Since this is a limited study target, 
there are still other convex pentagons with a finite Heesch number. Indeed, we 
have also found some convex pentagons with $H(T) = 1$, whose corona cannot 
be edge-to-edge. (See Appendix. Contrary to our expectations, many convex
 pentagons with a finite Heesch number were found. Therefore, in our search, 
we decided to temporarily exclude convex pentagons with a finite Heesch 
number whose first corona cannot be edge-to-edge. Therefore, we searched 
only the convex pentagons with finite Heesch number whose first corona 
can be edge-to-edge).

All convex pentagons with a finite Heesch number that we have discovered are 
convex pentagons with $H(T) = 1$, but it is unknown whether convex pentagons 
(or convex polygons) with $2 \le H(T) < \infty $ exist.

\begin{figure}[htbp]
 \centering\includegraphics[width=14.5cm,clip]{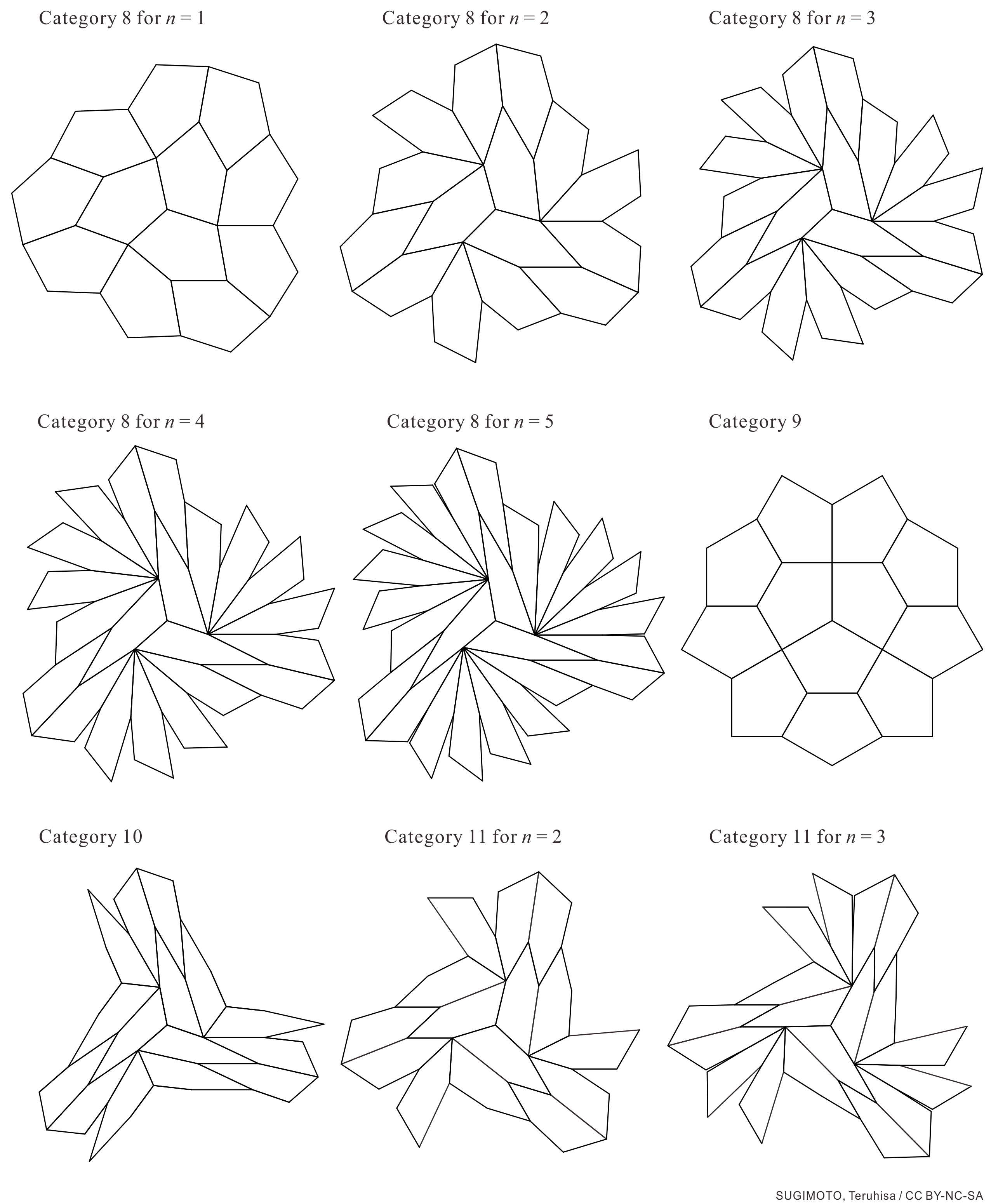} 
  \caption{{\small 
Cluster of three convex pentagons that are surrounded only one 
time using congruent copies of pentagons.
} 
\label{fig25}
}
\end{figure}

Finally we would like to present a new question based on this research 
result. We noticed that in the cases of Categories 8--11, the cluster of 
three convex pentagons can be surrounded only one time. Figure~\ref{fig25} 
shows these examples. Thus, we present the following question.

\begin{ques*} 
When a tile $T$ is convex and can be surrounded only one 
time using congruent copies of $T$ (i.e., $T$ is a convex tile with the Heesch 
number $1$), there is $T$ such that the cluster of three tiles can be surrounded 
only one time using congruent copies of $T$. For an integer $k \ge 4$, is there a 
convex tile $T$ with the Heesch number $1$ whose cluster of $k$ tiles can be 
surrounded only one time using congruent copies of $T$? Also, what is the 
maximum value of $k$?
\end{ques*}
 
\noindent
Note that Mann declared that for each integer $j \ge 1$, there is a tile 
$T_{j}$ such that two congruent copies of $T_{j}$ can surround $j$ copies of 
$T_{j}$ (i.e., the cluster of the $j$ tiles)~\cite{Mann_2002}. The tile $T_{j}$ was 
made by modifying the Voderberg tiles~\cite{Brass_2005, G_and_S_1987}, and 
is a concave tile with $H(T) = 1$.

Kaplan showed that there are infinite families of convex polygons with five 
or more edges with $H(T) = 1$. These convex polygons are shaped like ice 
cream cones or as such, bisected by lines of symmetry of polygons like an 
ice cream cone~\cite{Kaplan_2017_p3}. However, all of the convex polygons 
with $H(T) = 1$ that Kaplan showed in \cite{Kaplan_2017_p3} can form 
only non-edge-to-edge corona. These convex polygons obtained by 
bisecting a convex polygon such as an ice cream cone correspond to 
convex tiles with $H(T) = 1$, whose cluster of two tiles can be 
surrounded only one time using congruent copies. Therefore, there are 
infinite families of convex tiles with $H(T) = 1$ whose cluster of two tiles 
can be surrounded only one time. Of course, convex tiles with $H(T) = 1$, 
whose cluster of three tiles can be surrounded only one time using congruent 
copies, will be able to adjust state such that the cluster of two tiles can 
be surrounded only one time.

Our question above may be more interesting by limiting the shape of the tile 
to asymmetry, or by limiting the corona to edge-to-edge.

\renewcommand{\figurename}{{\small Figure.}}

\bigskip

\bigskip
\noindent
\textbf{Acknowledgments.} 
The authors thank Professor Shigeki Akiyama of University of Tsukuba for 
giving a chance to tackle the problem of convex pentagons with finite Heesch 
number. The authors would like to thank Professor Yoshio Agaoka of Hiroshima 
University and Professor Craig S. Kaplan of University of Waterloo for 
providing valuable comments.

\appendix
\def\thesection{Appendix }
\section{}

\begin{figure}[htbp]
 \centering\includegraphics[width=14.5cm,clip]{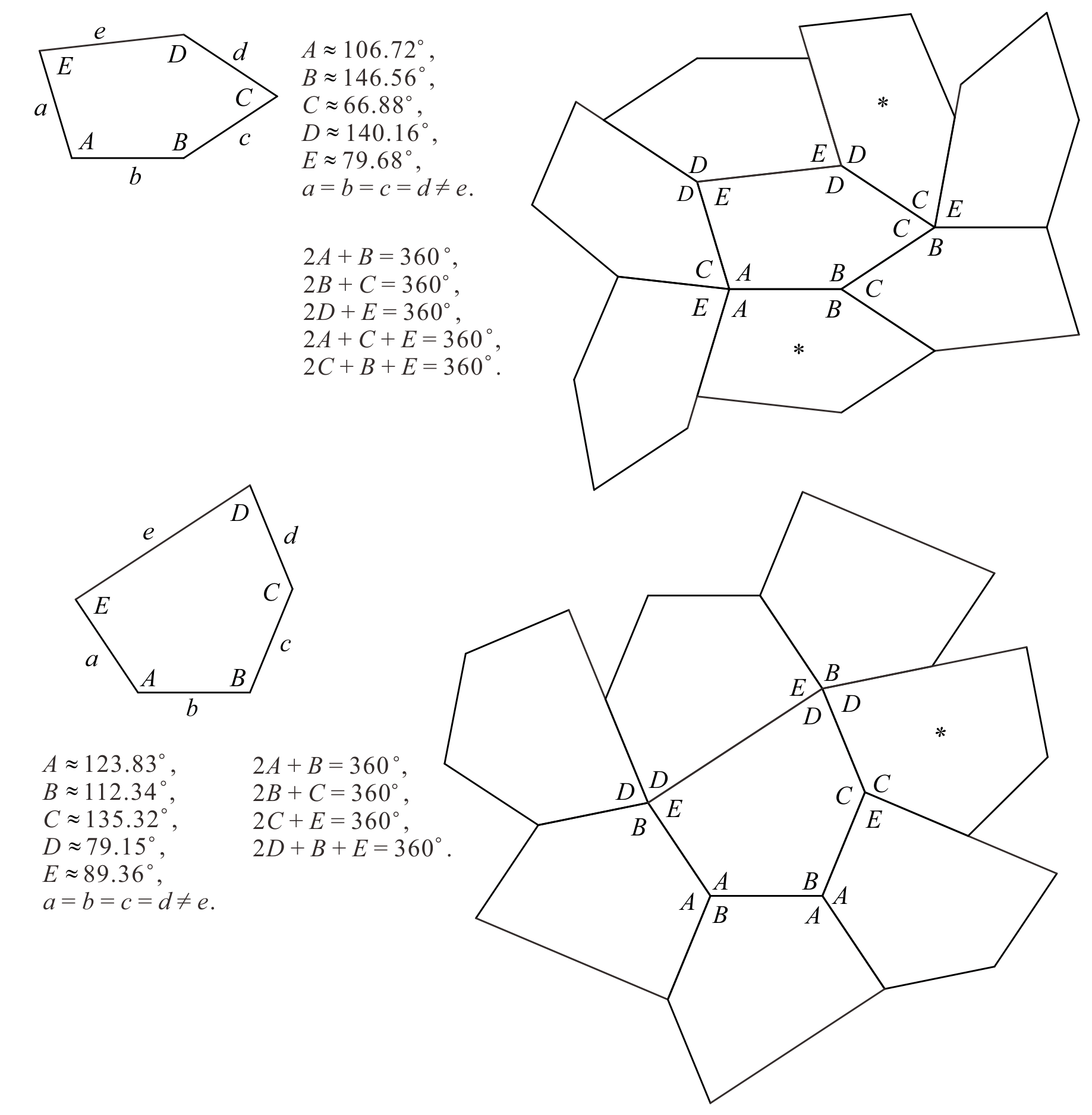} 
  \caption{{\small 
Convex pentagons with $H(T) = 1$ whose corona cannot be 
edge-to-edge.
} 
\label{fig26}
}
\end{figure}

Figure~\ref{fig26} shows an example of convex pentagons with $H(T) = 1$ whose corona 
cannot be edge-to-edge that we discovered.

\renewcommand{\figurename}{{\small Figure.}}

\end{document}